\documentclass{amsart}
\usepackage{amsmath,amsthm}%,amsfonts}
\usepackage {latexsym}
\usepackage{amssymb}

\newcommand{\beal}{\begin{align}}
\newcommand{\enal}{\end{align}}
\newcommand{\bealn}{\begin{align*}}
\newcommand{\enaln}{\end{align*}}
\newcommand{\bear}{\begin{eqnarray}}
\newcommand{\eear}{\end{eqnarray}}
\newcommand{\beeq}{\begin{equation}}
\newcommand{\eneq}{\end{equation}}

\newcommand{\eps}{{\varepsilon}}
\newcommand{\R}{{\mathbb R}}

\newcommand{\tileps}{{\tilde{\eps}}}

\def\bm{\left[ \begin{array}{cc}}
\def\endm{\end{array}\right]}

\def\eps{\varepsilon}

\def\bm{\left[\begin{matrix} }
\def\endm{\end{matrix}\right]}

\def\R{{\mathbb R}}

\newtheorem{theorem}{Theorem}
\newtheorem{lemma}[theorem]{Lemma}
\newtheorem{defi}[theorem]{Definition}
\newtheorem{cor}[theorem]{Corollary}
\newtheorem{prop}[theorem]{Proposition}

\theoremstyle{remark}
\newtheorem{remark}[theorem]{Remark}

\renewcommand{\Im}{\,{\rm Im}\,}
\renewcommand{\Re}{\,{\rm Re}\,}

\renewcommand{\hat}{\widehat}
\renewcommand{\epsilon}{\eps}
\renewcommand{\tilde}{\widetilde}
\numberwithin{equation}{section}
\numberwithin{theorem}{section}

\begin{document}

\title[wave maps full blow up range]{Optimal polynomial blow up range for critical wave maps.}

\author{Can Gao, Joachim Krieger}



\thanks{Support of the Swiss National Fund for
the second author is gratefully acknowledged.}

\begin{abstract}
We prove that the critical Wave Maps equation with target $S^2$ and origin $\R^{2+1}$ admits energy class blow up solutions of  the form 
\[
u(t, r) = Q(\lambda(t)r) + \eps(t, r)
\]
where $Q:\R^2\rightarrow S^2$ is the ground state harmonic map and $\lambda(t) = t^{-1-\nu}$ for any $\nu>0$. This extends the work \cite{KST0}, where such solutions were constructed under the assumption $\nu>\frac{1}{2}$. In light of a result of Struwe \cite{Struwe1}, our result is optimal for polynomial blow up rates.  
\end{abstract}

\maketitle

\section{Introduction}

This paper considers the issue of obtaining the optimal polynomial range of blow up dynamics for critical co-rotational Wave Maps from $\R^{2+1}$ into $S^2$, the standard two-dimensional sphere. 
Recall that a map 
\[
u: \R^{2+1}\longrightarrow S^2
\]
is considered a {\it{Wave Map}}, provided it is formally critical with respect to the (formal) Lagrangian action functional 
\[
\mathcal{L}(u): = \int_{\R^{2+1}}\langle\partial_{\alpha}u, \partial^{\alpha}u\rangle_{S^2}\,d\sigma,\,\partial^{\alpha} = m^{\alpha\beta}\partial_{\beta}
\]
 where, with $\alpha = 0,1,2$ space-time indices, the Einstein summation convention is in force, and $m^{\alpha\beta}$ is the Minkowski metric with signature $(-1,1,1)$. 
 Wave Maps from a $2+1$-dimensional background are {\it{energy critical}}, meaning that the natural conserved energy 
 \begin{equation}\label{eq:energy}
 \mathcal{E}(u): = \int_{\R^2}\big[|u_t|^2 + |\nabla_x u|^2|\big]\,dx
 \end{equation}
 is invariant under the intrinsic scaling 
 \[
 u(t, x)\rightarrow u(\lambda t, \lambda x)
 \]
 for the Wave Maps equation. 
 The Wave Maps equation has a remarkable so-called null-structure, as evidenced by its explicit form 
 \begin{equation}\label{eq:WMS2}
 \Box u = -u_{tt} + \triangle u = -u(-|u_t|^2 + |\nabla_x u|^2),\,u(t, x)\in S^2\subset \R^3
 \end{equation}
 This null-structure is responsible for the fact that \eqref{eq:WMS2} enjoys an almost optimal local well-posedness property: from \cite{KlMa}, it is known that \eqref{eq:WMS2} is strongly locally well-posed (in the sense of real analytic dependence of the solution on the data) in any space $H^{s}$, $s>1$.
On the other hand, from \cite{dAncGeo}, it is known that \eqref{eq:WMS2} is {\it{ill-posed}} (however, only in the sense of non-uniform continuous dependence of a local solution on the data) in any $H^{s}, s<1$. 
In the delicate borderline case of data in $H^1$ (corresponding to the energy \eqref{eq:energy}
), it is known\footnote{For an earlier result in the equivariant context, see \cite{ShaTah}.}, see \cite{T2}, and more recently \cite{SterbTat1}, that for $s>1$, {\it{ $H^s$-smooth data of small enough energy result in a global $H^s$-smooth solution}}. Furthermore, the solutions scatter at infinity like free waves, provided the initial data are $C^\infty$-smooth and  constant outside of a compact set, say. 
 In fact, the recent result \cite{SterbTat1} furnishes the optimal energy threshold, namely that of the minimum energy non-trivial harmonic map $Q$ from $\R^2\rightarrow S^2$, without any symmetry assumptions on the map. An earlier result \cite{CoKeMe} derived such a result in the co-rotational context. See also \cite{CoKeLaSch1},  \cite{CoKeLaSch2} for developments in the context of energy much above the ground state.  
Since the work \cite{KST0}, and later \cite{RaRod}, it has been known that for any $\eps>0$, there exist initial data\footnote{They may be chosen of any regularity $H^s, s>1$.} of energy $\mathcal{E}(Q) + \eps$ and which lead to finite time singularity formation. See also \cite{RodSterb} for blow up solutions with energy $>4\mathcal{E}(Q)$.  
In fact, the works \cite{KST0}, \cite{RaRod}, produced different blow up rates, the former exhibiting a continuum of blow up rates, the latter a more rigid rate but in turn demonstrably stable (within the co-rotational class). 
To explain this further, we recall the fundamental work \cite{Struwe1} by M. Struwe on the structure of singularities. Struwe shows that if 
\[
u: [0, T)\times \R^{2}\longrightarrow S^2
\]
is a smooth co-rotational\footnote{This means that if one uses spherical coordinates on $S^2$, and polar coordinates on the plane $\R^2$ of spatial variables, then the wave map can be described by $(t, r, \theta)\longrightarrow (\theta, u(t, r))$.} wave map which cannot be smoothly extended past time $T$, then there exists a sequence of times  $t_i\rightarrow T$  as well as a sequence of parameters $\lambda_i\rightarrow + \infty$ with the property that on each fixed time slice $t = t_i$, we can write
\[
u(t_i, x) = Q(\lambda(t_i)x) + \eps(t_i, x)
\]
 where $Q$ represents the ground state co-rotational harmonic map $Q: \R^2\rightarrow S^2$, while the error term $\eps$ satisfies 
 \[
\lim_{i\rightarrow\infty}\mathcal{E}_{loc}\big(\eps(t_i, x)\big): = \lim_{i\rightarrow\infty}\int_{|x|<t_i}\big[|\partial_t\eps(t_i, x)|^2 +|\nabla_x\eps(t_i, x)|^2\big]\,dx = 0
 \]
 Furthermore, Struwe established an upper bound on the blow up rate 
 \[
 \lim_{i\rightarrow\infty}\lambda(t_i)(T-t_i) = +\infty
 \]
The blow up rates exhibited in \cite{KST0}, \cite{RaRod}, of course obey this asymptotic, and in fact we have 
\[
\lambda(t) = (T-t)^{-\nu-1}
\]
with $\nu >\frac{1}{2}$ for the solutions constructed in \cite{KST0}. It then remains a very natural question to decide whether in fact all $\nu>0$ are  admissible. In this paper, we provide a positive answer to this.
To formulate the main theorem, we recall that co-rotational wave maps may be parametrized in terms of a function $u(t, r)\rightarrow \R$ which solves the {\it{scalar wave equation}} 
\begin{equation}\label{eq:WMcorot}
-\partial_{tt}u + \partial_{rr}u + \frac{1}{r}\partial_r u = \frac{\sin(2u)}{2r^2}
\end{equation}
In terms of this representation, the ground state harmonic map (which corresponds to a static wave map) is given by 
\[
Q(r) = 2\arctan r
\]
The function $u(t, r)$ is to be thought of as a function on $\R^2$, thus the conserved energy is given by 
\[
\int_0^{\infty}\big[u_t^2 + |u_r|^2 + \frac{\sin^2(u)}{r^2}\big]\,rdr
\]

\begin{theorem}\label{thm:Main} For any $\nu>0$, there exist $T>0$ and co-rotational initial data $(f, g)$ with 
\[
(f-\pi, g)\in H_{\R^2}^{1+\frac{\nu}{2}-}\times H_{\R^2}^{\frac{\nu}{2}-}
\]
which result in a\footnote{Here we use the identification of the wave map with a function $u(t, r)$ as before.} solution $u(t, r)$, $t\in (0, T]$ which blows up at time $t = 0$ and has the following representation: 
\[
u(t, r) = Q(\lambda(t)r) + \eps(t, r)
\]
where $\lambda(t) = t^{-1-\nu}$, and such that the function 
\[
(\theta, r) \longrightarrow \big(e^{i\theta}\eps(t, r), e^{i\theta}\eps_t(t, r)\big)\in H^{1+\nu-}(\R^2)\times H^{\nu-}(\R^2)
\]
uniformly in $t$. Also, we have the asymptotic as $t\rightarrow 0$
\[
\mathcal{E}_{loc}\big(\eps(t, \cdot)\big)\lesssim (t\lambda(t))^{-1}\log^2 t
\]
\end{theorem} 

\section{Some remarks on the result}

Our approach to the theorem is following closely the one in \cite{KST0}, with a key modification in the second part which essentially follows \cite{KS}. 
Specifically, we recall that the construction in \cite{KST0} has two essentially distinct stages: 
\begin{itemize}
\item In a first stage, we construct an approximate solution, denoted by 
\[
u_{approx}(t, r) = Q(\lambda(t)r) + u^e(t, r)
\]
where the correction term $u^e(t, r)$ is obtained by iteratively solving certain 'elliptic approximations' to the wave equation \eqref{eq:WMcorot}. While $u_{approx}(t, r)$ is not an exact solution of \eqref{eq:WMcorot}, it is a very accurate solution, in that we can ensure that given $N\geq 0 $, we can ensure that the error 
\[
-\partial_{tt}u_{approx} + \partial_{rr}u_{approx} + \frac{1}{r}\partial_r u_{approx} - \frac{\sin(2u_{approx})}{2r^2} = O(t^N).
\]
Of course  the larger $N$, the more 'elliptic correction terms' need to be added to $Q(\lambda(t)r)$.  It is important to observe here that the restriction $\nu>\frac{1}{2}$ imposed in \cite{KST0} does not come in at this stage; in fact, any $\nu>0$ will suffice. 
\item In a second stage, we complete the approximate solution $u_{approx}$ to an exact one by adding a correction term $\eps(t, r)$. This latter correction term is now determined by solving an actual wave equation, albeit one with a time dependent potential term. Dealing with the latter forces one to develop some rather sophisticated spectral theory. 
To find $\eps$, one implements a fixed point argument in a suitable Banach space, and it is here, in the treatment of the nonlinear terms with singular weights, that the restriction on $\nu$ comes in. Indeed, in Lemma 8.5 in \cite{KST0}, the bound (notation to be explained further below) 
\[
\|R^{-\frac{3}{2}}fg\|_{H_{\rho}^{\alpha+\frac{1}{4}}}\lesssim \|f\|_{H_{\rho}^{\alpha+\frac{1}{2}}}\|g\|_{H_{\rho}^{\alpha+\frac{1}{2}}}
\]
is derived which holds provided $\alpha>\frac{1}{4}$. Since the iterates for $\eps$ live naturally in the space $H_{\rho}^{\frac{1}{2}+\frac{\nu}{2}-}$, the condition $\nu>\frac{1}{2}$ used in \cite{KST0} follows. 
\end{itemize}

In the present work, we overcome this restriction as follows: 
\begin{itemize}
\item First, we analyze the 'zeroth iterate' (to be explained below) for (a suitable variant of)$\eps$, and show that we can split this into the sum of two terms, one of which has a regularity gain which lands us in the regime where the Lemma 8.5 in \cite{KST0} is applicable, the other of which does not gain regularity but satisfies an a priori $L^\infty$-bound near the symmetry axis $R = 0$. Note that the regularity requirement in Lemma 8.5 in \cite{KST0} comes primarily from the singular weight $R^{-\frac{3}{2}}$ at $R = 0$, and so an a priori bound on the (weighted) $L^\infty$ norm will be seen to suffice to estimate an expression such as $R^{-\frac{3}{2}}\eps^2$. Intuitively, the reason why we can control the part of the zeroth iterate near $R = 0$ comes from the fact that the singular behavior of the approximate solution from the first part of the construction and the error it generates is localized to the boundary of the light cone. 
\item Second, by writing the equation for the distorted Fourier transform of (a variant of ) $\eps$ in a way that subtly differs from the one in \cite{KST0}, we manage to show that the higher iterates all differ from the zeroth iterate by terms with a smoothness gain. This will then suffice to show the desired convergence. 
\end{itemize} 

\section{Construction of an approximate solution}\label{sec:approxsol}

Here we shall follow closely the procedure in \cite{KST0}, but also correct for certain (inessential) algebraic errors in the latter reference. In particular, we shall slightly modify the function spaces used (again without any major consequence). Denote 
\[
R = \lambda(t)r,\,\lambda(t) = t^{-1-\nu},\,\nu>0
\]
Also, write $u_0(R): = Q(R) = 2\arctan R$. 
We state the following, quite analogous to the result in \cite{KST0}:  
\begin{theorem}\label{thm:approxsol} Assume $k\in \mathbf{N}$. There exists an approximate solution $u_{2k-1}(R)$ for \eqref{eq:WMcorot} which can be written as 
\[
u_{2k-1}(t, r) = Q(R) + \frac{c_k}{(t\lambda)^2}R\log(1+R^2) +\frac{\tilde{c}_k}{(t\lambda)^2}R  + O\big(\frac{(\log(1+R^2))^2}{(t\lambda)^2}\big)
\]
with a corresponding error of size\footnote{The extra factor $(1-\frac{R}{\lambda t})^{-\frac{1}{2}}$ here arises for $\nu<\frac{1}{2}$, and is not present in \cite{KST0}.} 
\[
e_{2k-1} = (1-\frac{R}{\lambda t})^{-\frac{1}{2}+\nu}\big(\frac{R(\log(1+R^2))^2}{(t\lambda)^{2k}}\big)
\]
Here the implied constant in the $O(\ldots)$ symbols are uniform in $t\in (0,\delta]$ for some $\delta = \delta(k)>0$ sufficiently small. 
\end{theorem}
The construction of this solution follows very closely the treatment in \cite{KST0}. Specifically, we shall arrive at the $k$-th approximation by adding $k$ correction terms to $u_0$: 
\[
u_k = u_0 + \sum_{j=1}^k v_j
\]
Write 
\[
e_k =  \partial_t^2 u_k -\partial_r^2 u_k - \frac{1}{r}\partial_r u_k + \frac{\sin(2u_k)}{2r^2}
\]
From \cite{KST0} we recall how the correction terms $v_k$ are computed inductively: for each $k$, we employ a splitting 
\[
e_k = e^0_k + e^1_k
\]
where $e^1_k$ denotes certain higher order error terms relegated to a later stage of the inductive process. Then depending on whether $k$ is even or not, we define 
\begin{equation}\label{eq:odditerate}
\big(\partial_r^2 + \frac{1}{r}\partial_r - \frac{\cos(2u_0)}{r^2}\big)v_{2k+1} = e^0_{2k}
\end{equation}
\begin{equation}\label{eq:eveniterate}
\big(-\partial_t^2 + \partial_r^2 + \frac{1}{r}\partial_r - \frac{1}{r^2}\big)v_{2k} = e^0_{2k-1}
\end{equation}
where we impose trivial Cauchy data at $r = 0$,
resulting in the new error terms 
\[
e_{2k+1} = e^1_{2k} - \partial_t^2v_{2k+1} + N_{2k+1}(v_{2k+1}),\,e_{2k} = e^1_{2k-1} + N_{2k}(v_{2k})
\]
Here we have introduced the expressions 
\begin{align}\label{NLT}
&N_{2k}(v)=\frac{1-\cos(2u_{2k-1})}{r^2}v+\frac{\sin(2u_{2k-1})}{2r^2}(1-\cos(2v))+\frac{\cos(2u_{2k-1})}{2r^2}(2v-\sin(2v))\\
&N_{2k+1}(v)=\frac{\cos(2u_0)-\cos(u_{2k})}{r^2}v+\frac{\sin(2u_{2k})}{2r^2}(1-\cos(2v))+\frac{\cos(2u_{2k})}{2r^2}(2v-\sin(2v))
\end{align}

The key fact for this construction is that while \eqref{eq:eveniterate} is a wave equation, the ansatz that we will use to construct $v_{2k}$ will allow us to reformulate this problem as a singular elliptic Sturm-Liouville problem, which can be solved by standard ODE methods. It will then be seen that the errors are in fact decreasing near $t = 0$. The main challenge is to control the (increasingly complicated) corrections $v_k$ by placing them in suitable function spaces. 

We now define these spaces, implementing very subtle changes compared to \cite{KST0}, in the definition of the ingredients of $S^m(R^k(\log R)^l, \mathcal{Q}_n)$ below:

\begin{defi}\label{defi:Q_n}For $i\in \mathbf{N}$, let $j(i) = i$ if $\nu$ is irrational, respectively $j(i) = 2i^2$ if $\nu$ is rational. 
Then
\begin{itemize}
\item $\mathcal{Q}$ is the algebra of continuous functions $q: [0,1]\rightarrow \R$ with the following properties: 
\\
(i) $q$ is analytic in $[0,1)$ with even expansion around $a = 0$.\\
(ii) near $a = 1$ we have an absolutely convergent expansion of the form
\begin{align*}
q(a) = &q_0(a) + \sum_{i=1}^\infty (1-a)^{\beta(i)+\frac{1}{2}}\sum_{j=0}^{j(i)}q_{i,j}(a)\big(\log(1-a)\big)^j\\
&+ \sum_{i=1}^\infty (1-a)^{\tilde{\beta}(i)+\frac{1}{2}}\sum_{j=0}^{j(i)}\tilde{q}_{i,j}(a)\big(\log(1-a)\big)^j
\end{align*}
with analytic coefficients $q_0, q_{i,j}$, and $\beta(i) = i\nu$, $\tilde{\beta}(i) = \nu i+\frac{1}{2}$.  
\item $\mathcal{Q}_n$ is the algebra which is defined similarly, but also requiring $q_{i,j}(1) = 0$ if $i\geq 2n+1$. 
\end{itemize}
\end{defi}
We also define the space of functions obtained by differentiating $\mathcal{Q}_n$: 
\begin{defi}\label{defi:Q'_n}
Define $\mathcal{Q}'$ as in the preceding definition but replacing $\beta(i)$ by $\beta'(i): = \beta(i) - 1$, and similarly for $\mathcal{Q}'_n$.  
\end{defi}

The next definition also diverges slightly from the one in \cite{KST0}, see also \cite{KS}: 
\begin{defi}\label{defi:S^n} $S^n(R^k(\log R)^l)$ is the class of analytic functions $v: [0,\infty)\rightarrow \R$ with the following properties:\\
(i) $v$ vanishes of order $n$ at $R = 0$.\\
(ii) $v$ has a convergent expansion near $R = \infty$
\[
v = \sum_{\substack{0\leq j\leq l+i\\ i\geq 0}}c_{ij}R^{k-i}(\log R)^j
\]
\end{defi}
Next, introduce the symbols 
\[
b_1 = \frac{\big(\log(1+R^2)\big)^2}{(t\lambda)^2},\,b_2 = \frac{1}{(t\lambda)^2}
\]
The final function space is also slightly different than the one in \cite{KST0}: 
\begin{defi}\label{defi:S^mQ_n}  Pick $t$ sufficiently small such that both $b_1, b_2$, when restricted to the light cone $r\leq t$ are of size at most $b_0$. 
\begin{itemize}
\item $S^m(R^k(\log R)^l, \mathcal{Q}_n)$ is the class of analytic functions $v: [0,\infty)\times [0,1)\times [0,b_0]^2\rightarrow\R$ so that\\
(i) $v$ is analytic as a function of $R, b_1, b_2$, 
\[
v: [0,\infty)\times [0, b_0]^2\rightarrow {\mathcal{Q}}_n
\]
(ii) $v$ vanishes to order $m$ at $R = 0$.\\
(iii) $v$ admits a convergent expansion at $R = \infty$, 
\[
v(R,\cdot,b_1, b_2) = \sum_{\substack{0\leq j\leq l+i\\ i\geq 0}}c_{ij}(\cdot, b_1, b_2)R^{k-i}(\log R)^j  
\]
where the coefficients $c_{ij}: [0, b_0]^2\rightarrow \mathcal{Q}_n$ are analytic with respect to $b_{1,2}$. 
\item $IS^m(R^k(\log R)^l, \mathcal{Q}_n)$ is the class of analytic functions $w$ inside the cone $r<t$ which can be represented as 
\[
w(t, r) = v(R, a, b_1, b_2),\,v\in S^m(R^k(\log R)^l, \mathcal{Q}_n)
\]
and $t>0$ sufficiently small.  
\end{itemize}
\end{defi}

In the sequel, we shall show inductively that one can choose the corrections $v_k$ to satisfy the following: 
\begin{equation}\label{eq:v_2k-1}
v_{2k-1}\in \frac{1}{(t\lambda)^{2k}}IS^3\big(R(\log R)^{2k-1}, \mathcal{Q}_{k-1}\big)
\end{equation}
\begin{equation}\label{eq:e_2k-1}
t^2e_{2k-1}\in \frac{1}{(t\lambda)^{2k}}IS^1\big(R(\log R)^{2k-1}, \mathcal{Q}'_{k-1}\big)
\end{equation}
\begin{equation}\label{eq:v_2k}
v_{2k}\in \frac{1}{(t\lambda)^{2k+2}}IS^3\big(R^3(\log R)^{2k-1}, \mathcal{Q}_{k}\big)
\end{equation}
\begin{equation}\label{eq:e_2k}
t^2e_{2k}\in \frac{1}{(t\lambda)^{2k}}\big[IS^1\big(R^{-1}(\log R)^{2k}, \mathcal{Q}_{k}\big) + \langle b_1, b_2\rangle[IS^1\big(R(\log R)^{2k-1}, \mathcal{Q}'_{k}\big)
\end{equation}
and the starting error $e_0$ satisfying 
\[
e_0\in IS^1(R^{-1})
\]
Here we denote by $\langle b_1, b_2\rangle$ the ideal generated by $b_1, b_2$ inside the algebra generated by $b_1, b_2$. 
We first explicitly compute the first and second corrections $v_{1,2}$, and then automate the process for the higher iterates. To begin with, from the calculation in \cite{KST0}, we find 
\[
e_0 = \frac{1}{t^2}\big( (\nu+1)^2\frac{4R}{(1+R^2)^2} - \nu(\nu+1)\frac{2R}{1+R^2}\big)
\]

\subsection{The first correction}

If we try to make $u_1 = u_0 + \eps$ an exact solution, then $\eps$ needs to solve 
\begin{equation}\label{eq:epse_0}
\big(-\partial_{tt}+\partial_{rr}+\frac{1}{r}\partial_r\big)\varepsilon-\frac{\cos(2u_{0})}{2r^2}\sin{2\varepsilon}+\frac{\sin(2u_0)}{2r^2}(1-\cos(2\varepsilon))= e_0
\end{equation}
Introduce the operator 
\[
\tilde{\mathcal{L}}: = \partial_R^2 + \frac{1}{R}\partial_R - \frac{\cos(2u_0)}{R^2} = \partial_R^2 + \frac{1}{R}\partial_R - \frac{1}{R^2}\frac{1-6R^2 + R^4}{(1+R^2)^2}
\]
Now if we neglect the time derivatives $-\partial_{tt}\eps$ as well as the nonlinear term $\frac{\sin(2u_0)}{2r^2}(1-\cos(2\varepsilon))$ in \eqref{eq:epse_0} and replace the exact correction  $\eps$ by an approximate one $v_1$, we obtain the following relation 
\[
(t\lambda)^2\tilde{\mathcal{L}}v_1 = t^2 e_0
\]
which is a non-degenerate second order ODE and hence solvable by standard methods. 
Introduce the conjugated operator $\tilde{\mathcal{L}}$ by means of 
\[
-\mathcal{L}(\sqrt{R}v) = \sqrt{R}\tilde{\mathcal{L}}v
\]
Then one has 
\[
-\mathcal{L} = \partial_R^2 - \frac{3}{4R^2} + \frac{8}{(1+R^2)^2},
\]
and a fundamental system for the operator $\mathcal{L}$ is given by (see \cite{KST0})
\[
\phi(R)=\frac{R^{3/2}}{1+R^2}\quad \theta(R)=\frac{-1+4R^2\log R+R^4}{\sqrt{R}(1+R^2)}.
\]
With this choice, we have $W(\phi, \theta) = 2$. We have the variation of constants formula
\begin{align*}
(t\lambda)^2v_1 = \frac{1}{2}R^{-\frac{1}{2}}\theta(R)\int_0^R \phi(R')\sqrt{R'}f(R')\,dR' - \frac{1}{2}R^{-\frac{1}{2}}\phi(R)\int_0^R\theta(R')\sqrt{R'}f(R')\,dR'
\end{align*}
where we have put $f = t^2e_0$. Then compute for large $R$ and suitable constants $c_1, c_2, c_3, c_4, d_1, d_2, d_3, d_4$ 
\begin{align}
&R^{-1/2}\theta(R)\int_0^R\phi(R')\sqrt{R'}t^2e_0(R')dR'\\
&=\frac{-1+4R^2\log R+R^4}{R(1+R^2)}\int_0^R(c_1+\frac{c_2}{1+R'^2}) \Big(\frac{c_3}{1+R'^2}+\frac{c_4}{(1+R'^2)^2}\Big)d(1+R'^2)\nonumber\\
&=\frac{-1+4R^2\log R+R^4}{R(1+R^2)}\Big(\frac{c_1}{1+R^2}+\frac{c_2}{(1+R^2)^2}+c_3\log(1+R^2)+c_4\Big)\nonumber\\
&=d_1R\log R+d_2 R+d_3 R^{-1}\log^2 R + d_4 R^{-1} + O(R^{-2}\log^2 R). \nonumber
\end{align}
and similarly (with re-labelled coefficients)
\begin{align}
&R^{-1/2}\phi(R)\int_0^R\theta(R')\sqrt{R'}t^2e_0(R')dR'\\
&=\frac{R}{1+R^2}\int_0^R\frac{R'^4+4R'^2\log
  R'-1}{1+R'^2}\Big((\nu+1)^2\frac{4R'^3}{(1+R'^2)^2}- \nu(\nu+1)\frac{2R'}{1+R'^2}\Big)dR'\nonumber\\
&=\frac{R}{(1+R^2)}\int_{0}^R\Big(c_1+c_2(1+R'^2)+\frac{c_3R'^2\log
R'}{1+R'^2}\Big)\Big(\frac{c_5}{1+R'^2}+\frac{c_6}{(1+R'^2)^2}\Big)d(1+R'^2)\nonumber\\
&=R\Big(\sum_{i=-3}^0d_i(1+R^2)^i +d_3\log
(1+R^2)+\frac{d_4\log
(1+R^2)}{1+R^2}+\frac{d_5\log
(1+R^2)}{(1+R^2)^2}+\frac{d_6(\log(1+R^2))^2}{1+R^2}\Big)\nonumber\\
&+O(R^{-3}\log^2 R)\nonumber\\
&=e_1R\log R+e_2R+e_3\log R+e_4+ O(R^{-1}\log^2R)\nonumber
\end{align}
Furthermore, since $e_0$ vanishes to first order at $R=0$, it follows that $v_1$ vanishes to third order at zero, Combining these observations, we find that indeed
\[
v_1\in \frac{1}{(t\lambda)^2}IS^3(R\log R, \mathcal{Q}_0)
\]
as required from \eqref{eq:v_2k-1}. 

\subsection{The error generated after the first correction}

Her we calculate $t^2e_1$. This is given by 
\begin{align}\label{e1}
t^2e_1&=-t^2(-\partial^2_t+\partial^2_r+\frac 1r \partial_r)(u_0+v_1)+t^2\frac{\sin(2u_0+2v_1)}{2r^2}\\
&=t^2\Big[\partial_{tt}v_1-\frac{\sin
  2u_0}{2r^2}(1-\cos(2v_1))-\frac{\cos
  (2u_0)}{2r^2}(2v_1-\sin(2v_1))\Big]\nonumber\\
&=t^2\partial_{tt}v_1-\frac{\sin
  2u_0}{2R^2}(t\lambda)^2(1-\cos(2v_1))-\frac{\cos
  (2u_0)}{2R^2}(t\lambda)^2(2v_1-\sin(2v_1))\nonumber
\end{align}
Then we use that for $l\geq 1$
\[
R^{-2}(t\lambda)^2v_1^{2l+1} \in \frac{1}{(t\lambda)^2}IS^3(R\log R, \mathcal{Q}_0),\,R^{-2}(t\lambda)^2v_1^{2l} \in \frac{1}{(t\lambda)^2}IS^3(\log^2 R, \mathcal{Q}_0),
\]
which in addition to the fact that $u_0$ admits an expansion in terms of inverse powers of $R$ near $R = +\infty$ leads to 
\[
t^2 e_1\in \frac{1}{(t\lambda)^2}IS^3(R\log R, \mathcal{Q}_0)\subset \frac{1}{(t\lambda)^2}IS^3(R\log R, \mathcal{Q}'_0),
\]
as required. 

\subsection{The second correction}

Now we intend to add a second correction $v_2$ in order to reduce the error $e_1$ from the first stage. More precisely, this time we reduce this error near the light cone. 
Write $t^2 e_1$ in terms of its expansion at $R = \infty$: 
\[
t^2 e_1 = \frac{1}{(t\lambda)^2}\big[c_1 R\log R + c_2R + c_3\log R + c_4 + O(R^{-1}\log^2R)\big]
\]
for suitable coefficients $c_1,\ldots, c_4$. 
Neglecting the higher order error terms $O(R^{-1}\log^2R)$, we have to solve the equation 
\[
t^2\big(-\partial_t^2 + \partial_r^2 + \frac{1}{r}\partial_r - \frac{1}{r^2}\big)v_2 = t^2 e_1^0,
\]
where we write 
\[
t^2e_1^0: = \frac{1}{(t\lambda)^2}\big[c_1 R\log R + c_2R + c_3\log R + c_4\big] 
\]
Homogeneity considerations suggest making the following ansatz: $v_2 = w_2 + \tilde{w}_2$, where 
\[
w_2=\frac{1}{t\lambda}(W_2^1(a)\log R+W_2^0(a)),\quad \tilde{w}_2=\frac{1}{(t\lambda)^2}(\tilde{W}_2^1(a)\log
R+\tilde{W}_2^0(a)).
\]
To obtain the equations for the functions $W_2^1(a)$, we match powers of $R$ and $\log R$.  
We arrive at the following equations: 
\begin{align}\label{u_2W}
&t^2\tilde{\Box}(\frac{1}{t\lambda}W_2^i(a))=\frac{1}{t\lambda}(ac_{i+1}-F_i(a)),\quad i=1,0\\
&t^2\tilde{\Box}(\frac{1}{(t\lambda)^2}\widetilde{W}_2^i(a))=\frac{1}{(t\lambda)^2}(c_{i+2}-\widetilde{F}_i(a)),\quad i=1,0
\end{align}
where
\[
\tilde{\Box} = -\partial_t^2 + \partial_r^2 + \frac{1}{r}\partial_r - \frac{1}{r^2}
\]
as well as 
\begin{align*}
&F_1(a)=0,\quad
F_0(a)=((\nu+1)\nu+a^{-2})W_2^1(a)+(a^{-1}-(1+\nu)a)\partial_aW_2^1(a)\\
&\widetilde{F}_1(a)=0,\quad
\widetilde{F}_0=(2(\nu+1)\nu+a^{-2})\widetilde{W}_2^1(a)+(a^{-1}-(1+\nu)a)\partial_a\widetilde{W}_2^1(a)
\end{align*}
We conjugate out the power of $t$ and rewrite the equations in the $a$ variable
\begin{align*}
&\mathcal{L}_{\nu}W_2^i(a)=ac_{i+1}-F_i(a)\\
&\mathcal{L}_{2\nu}\widetilde{W}_2^i(a)=c_{i+2}-\widetilde{F}_i(a)
\end{align*}
where the one parameter family of operators $\mathcal{L}_{\beta}$ is
defined by
\begin{align}\label{dfnLb}
\mathcal{L}_{\beta}:=
(1-a^2)\partial^2_{\alpha}+(a^{-1}+2a\beta-2a)\partial_{a}+(-\beta^2+\beta
-a^{-2})
\end{align}
From \cite{KST0}, we know that there exist analytic solutions $W^i_{2}(a), \widetilde{W}_2^i(a)$ for \eqref{u_2W} on $[0,1)$, such that 
\[
W^i_{2}(a),\,i = 0,1,
\]
admits an odd power expansion around $a =0$ starting with $a^3$, while $\widetilde{W}_2^i(a)$ admits an even expansion around $a =0$, starting with $a^2$. 
Moreover, for $a$ near $1$, as shown in \cite{KST0}, we have expansions 
\[
W_2^1(a) = g_0(a) + g_1(a)(1-a)^{\nu+\frac{1}{2}} + g_2(a)(1-a)^{\nu+\frac{1}{2}}\log(1-a)
\]
\[
W_2^0(a) = h_0(a) + (1-a)^{\nu+\frac{1}{2}}\sum_{l=0}^2 h_{l+1}(a)[\log(1-a)]^l + (1-a)^{2\nu+1}h_{l+4}(a)[\log(1-a)]^l,
\]
where we have taken into account the most general case (when $\nu$ is irrational, there are fewer terms in the expansion). The result for $\tilde{W}_2^{1.0}(a)$ is of course analogous.
The expressions for $w_2, \tilde{w}_2$ are not quite what we want, since we need ultimately functions which vanish to odd order at $R = 0$, in order to ensure the desired smoothness. Furthermore, we also have the logarithmic factors $\log R$, which of course become singular at $R = 0$. 
In order to deal with these issues, we now {\it{re-define}} the correction terms $w_2, \tilde{w}_2$ in the following manner: 
\begin{align*}
&w_2=\frac{1}{t\lambda}(W_2^1(a)\frac{1}{2}\log (1+R^2)+W_2^0(a)),\\&\tilde{w}_2=\frac{1}{(t\lambda)^2}\frac{R}{(1+R^2)^{\frac{1}{2}}}(\tilde{W}_2^1(a)\frac{1}{2}\log(1+
R^2)+\tilde{W}_2^0(a)).
\end{align*}
Writing 
\[
\frac{1}{t\lambda}W_2^{1,0}(a) = \frac{1}{(t\lambda)^2}R Z_2^{1,0}(a)
\]
where now $Z_2^{1,0}(a)\in \mathcal{Q}_1$, while from construction we have $\tilde{W}_2^0(a)\in \mathcal{Q}_1$, and observing that $Z_2^{1,0}(a), \tilde{W}_2^0(a)$ vanish quadratically at $a = 0$ we see that 
\[
v_2 = w_2 + \tilde{w}_2\in \frac{1}{(t\lambda)^4}IS^3(R^3\log R, \mathcal{Q}_1),
\]
as required. 

\subsection{The error generated after the second correction $v_2$}

We write $u_2 = u_1 + v_2 = u_0 + v_1 + w_2 + \tilde{w}_2$, and need to estimate 
\begin{align*}
t^2e_2 = t^2(e_1 - e_1^0) + t^2\big(-\partial_t^2 + \partial_r^2 + \frac{1}{r}\partial_r - \frac{1}{r^2}\big)v_2 - t^2e_1^0 + t^2N_2(v_2)
\end{align*}
We check that each of the terms on the right satisfies the property \eqref{eq:e_2k} with $k = 1$. 
\\

{\it{(1): The contribution of $t^2(e_1 - e_1^0)$.}} From our choice of $e_1^0$, we immediately get 
\[
t^2(e_1 - e_1^0)\in \frac{1}{(t\lambda)^2}IS^1(R^{-1}(\log R)^2, \mathcal{Q}_1)
\]

{\it{(2): The contribution of $ t^2\big(-\partial_t^2 + \partial_r^2 + \frac{1}{r}\partial_r - \frac{1}{r^2}\big)v_2 - t^2e_1^0$.}}
This error is produced by replacing $\log R$ by $\frac{1}{2}\log(1+R^2)$, as well as by including the factor $\frac{R}{(1+R^2)^{\frac{1}{2}}}$. Thus we write this contribution as\footnote{Recall that $\tilde{\Box} = -\partial_t^2 + \partial_r^2 + \frac{1}{r}\partial_r - \frac{1}{r^2}$.} 
\begin{align*}
&t^2\tilde{\Box}\big[\frac{1}{t\lambda}W_2^1(a)(\frac{1}{2}\log (1+R^2) - \log R)\big]\\
&+t^2\tilde{\Box}\big[\frac{1}{(t\lambda)^2}\frac{R}{(1+R^2)^{\frac{1}{2}}}(\tilde{W}_2^1(a)(\frac{1}{2}\log(1+
R^2) - \log R)\big]\\
&+t^2\tilde{\Box}'\big[\frac{1}{(t\lambda)^2}\frac{R}{(1+R^2)^{\frac{1}{2}}}\tilde{W}_2^1(a)\frac{1}{2}\log(1+
R^2)\big]\\
&+\big(\frac{R}{(1+R^2)^{\frac{1}{2}}} - 1\big)t^2 e_1^0
\end{align*}
where the notation $\tilde{\Box}'$ means that at least one derivative falls on the factor $\frac{R}{(1+R^2)^{\frac{1}{2}}}$. 
Since $\frac{1}{2}\log (1+R^2) - \log R = O(R^{-2})$ as $R\rightarrow\infty$, and since $W_2^1(a)$ vanishes to third order at $a = 0$, we obtain easily that the first three expressions are in the space 
\[
\frac{1}{(t\lambda)^2}IS^1(R^{-1}, \mathcal{Q}_1')
\]
and since $\frac{R}{(1+R^2)^{\frac{1}{2}}} - 1 = O(R^{-2})$, the same applies to the last term above. 
This is not quite of the form required for \eqref{eq:e_2k}. However, we can rectify this by writing as in \cite{KST0} for any $t^2e\in \frac{1}{(t\lambda)^2}IS^1(R^{-1}, \mathcal{Q}_1')$
\[
t^2e = (1-a^2)t^2e + \frac{R^2}{(t\lambda)^2}t^2e
\]
This implies 
\[
\frac{1}{(t\lambda)^2}IS^1(R^{-1}, \mathcal{Q}_1')\subset \frac{1}{(t\lambda)^2}IS^1(R^{-1}, \mathcal{Q}_1) + b_2\frac{1}{(t\lambda)^2}IS^1(R, \mathcal{Q}_1')
\]
{\it{(3): The contribution of $t^2 N_{2}(v_2)$.}} Recall from \eqref{NLT} that we need to control three terms. First, we have  
\begin{align*}
&t^2\frac{1-\cos
  (2u_1)}{r^2}v_2\\&=\frac{(t\lambda)^2}{R^2}\Big(S^1(R^{-1}, \mathcal{Q}_1)+\frac{1}{(t\lambda)^2}S^3(R\log
R,  \mathcal{Q}_1)\Big)^2\times \frac{1}{(t\lambda)^4}S^3(R^3(\log R),  \mathcal{Q}_1)\\
&\in \frac{1}{(t\lambda)^2}\Big(S^3(R^{-1}(\log
R),  \mathcal{Q}_1)+\frac{1}{(t\lambda)^2}S^5(R(\log
R)^2,  \mathcal{Q}_1)+\frac{1}{(t\lambda)^4}S^7(R^3(\log R)^3,  \mathcal{Q}_1)\Big)\\
& \subset \frac{1}{(t\lambda)^2}\Big(S^3(R^{-1}(\log
  R),  \mathcal{Q}_1)+\frac{\langle b_1, b_2\rangle}{(t\lambda)^2}S^5(R(\log R),  \mathcal{Q}_1)\Big),
\end{align*}
as required. Further, just as in \cite{KST0}, one checks that  
\begin{align*}
t^2\frac{\sin (2u_1)}{2r^2}(1-\cos (2v_2))\in
\frac{1}{(t\lambda)^2}(S^1(R^{-1}(\log R)^2, \mathcal{Q}_1)+\langle b_1, b_2\rangle S^3(R(\log R),  \mathcal{Q}_1))
\end{align*}
and finally for the the cubic term
\begin{align*}
t^2\frac{\cos (2u_1)}{2r^2}(2v_2-\sin (2v_2))\in
\frac{\langle b_1, b_2\rangle}{(t\lambda)^2}S^1(R(\log R),  \mathcal{Q}_1).
\end{align*}
Combining all we have now, we conclude
$$t^2e_2\in \frac{1}{(t\lambda)^2}\big[S^1(R^{-1}(\log
R)^2, \mathcal{Q}_1)+\langle b_1, b_2\rangle S^1(R(\log R),  \mathcal{Q}_1)\big],$$
thus verifying \eqref{eq:e_2k} for $k = 1$. 

\subsection{The higher order corrections $v_k$, $k\geq 3$.} Here we repeat the preceding steps, making sure that the corrections and errors remain in the appropriate function spaces. We essentially use the same inductive procedure as in \cite{KST0}, with the same subtle changes as before. We do this in a number of steps:
\\

{\bf{Step 1}}: {\it{Given $u_{2k-2}$ with generating error $e_{2k-2}, k\geq 2$, as in \eqref{eq:e_2k},
choose $v_{2k-1}$ as in \eqref{eq:v_2k-1} with error $e_{2k-1}$ satisfying \eqref{eq:e_2k-1}.}} 
\\

This is accomplished exactly as in {\bf{Steps 1,2}} in \cite{KST0}. 
\\

{\bf{Step 2}}: {\it{Given $e_{2k-1}$ as in \eqref{eq:e_2k-1}, construct $v_{2k}$ as in \eqref{eq:v_2k}}}
\\

Here we have to diverge slightly from \cite{KST0}, since our definition of the algebra $S^m(R^l\log R^l)$ is different. Thus assume 
\[
t^2e_{2k-1}\in \frac{1}{(t\lambda)^{2k}}IS^1(R(\log R)^{2k-1},\,\mathcal{Q}'_{k-1})
\]
is given. We begin by isolating the leading component $e_{2k-1}^0$ which includes the terms of top degree in $R$ as well as those of one degree less (which is where we differ from \cite{KST0}). Thus we write 
$$t^2e^0_{2k-1}=\frac{R}{(t\lambda)^{2k}}\sum_{j=0}^{2k-1}q_j(a)(\log
R)^j+\frac{1}{(t\lambda)^{2k}}\sum_{j=0}^{2k}\tilde{q}_j(a)(\log
R)^j,\quad q_j,\tilde{q}_j\in \mathcal{Q}'_{k-1}$$
which we can rewrite as
$$t^2e^0_{2k-1}=\frac{1}{(t\lambda)^{2k-1}}\sum_{j=0}^{2k-1}aq_j(a)(\log
R)^j+\frac{1}{(t\lambda)^{2k}}\sum_{j=0}^{2k}\tilde{q}_j(a)(\log
R)^j$$
Consider the following equation
$$t^2\tilde{\Box}(v_{2k})=t^2e^0_{2k-1}.$$
Homogeneity considerations suggest that we should look for a solution
$v_{2k}$ which has the form
$$v_{2k}=\frac{1}{(t\lambda)^{2k-1}}\sum_{j=0}^{2k-1}W_{2k}^j(a)(\log R)^j+\frac{1}{(t\lambda)^{2k}}\sum_{j=0}^{2k}\widetilde{W}_{2k}^j(a)(\log R)^j$$
The one-dimensional equations for $W^j_{2k}$, $\widetilde{W}^j_{2k}$
are obtained by matching the powers of $\log R$. We get the following
systems
\begin{align*}
&t^2\tilde{\Box}\big( \frac{1}{(t\lambda)^{2k-1}}W_{2k}^j(a)\big)=\frac{1}{(t\lambda)^{2k-1}}(aq_j(a)-F_j(a))\\
&t^2\tilde{\Box}\big( \frac{1}{(t\lambda)^{2k}}\widetilde{W}_{2k}^i(a)\big)=\frac{1}{(t\lambda)^{2k}}(\tilde{q}_i(a)-\widetilde{F}_i(a))
\end{align*}
where $F_j(a)$, $\widetilde{F}_i(a)$ are
\begin{align*}
&F_j(a)=(j+1)[((\nu+1)\nu(2k-1)+a^{-2})W^{j+1}_{2k}+(a^{-1}-(1+\nu)a\partial_aW_{2k}^{j+1})]\\
&+(j+2)(j+1)((\nu+1)^2+a^{-2})W^{j+1}_{2k}\\
&\widetilde{F}_i(a)=(i+1)[(2(\nu+1)\nu
k+a^{-2})W^{i+1}_{2k}+(a^{-1}-(1+\nu)a\partial_aW_{2k}^{i+1})]\\
&+(i+2)(i+1)((\nu+1)^2+a^{-2})W^{i+1}_{2k}
\end{align*}
Here we make the convention that $W^j_{2k}(a), \widetilde{W}^i_{2k}=0$
for $j\geq 2k$ and $i\geq 2k+1$. Then we solve the systems
successively for decreasing values of $j,i$.
Conjugating out the power
of $t$ we get
\begin{align*}
&t^2\Big(-\Big(\partial_t+\frac{(2k-1)\nu}{t}\Big)^2+\partial^2_r+\frac
1r\partial_r-\frac{1}{r^2}\Big)W_{2k}^j(a)=aq_j-F_j(a)\\
&t^2\Big(-\Big(\partial_t+\frac{2k\nu}{t}\Big)^2+\partial^2_r+\frac
1r\partial_r-\frac{1}{r^2}\Big)\widetilde{W}_{2k}^i(a)=\tilde{q}_i-\widetilde{F}_i(a)
\end{align*}
with the definition of $\mathcal{L}_{\beta}$ we give in
$(\ref{dfnLb})$, we rewrite them in the $a$ variable
\begin{align*}
&\mathcal{L}_{(2k-1)\nu}W^j_{2k}=aq_j(a)-F_j(a)\\
&\mathcal{L}_{2k\nu}\widetilde{W}^i_{2k}=\tilde{q}_i(a)-\widetilde{F}_i(a)
\end{align*}
we claim that solving this system with Cauchy data at $a=0$ yields
solutions which satisfy
\begin{align*}
&W_{2k}^j(a)\in a^3\mathcal{Q}_{k},\quad j=\overline{0,2k-1}\\
&\widetilde{W}_{2k}^i\in a^2\mathcal{Q}_k,\quad i=\overline{0,2k}
\end{align*}
and this claim is established as in the computation of $v_2$ above,
see \cite{KST0}, lemma 3.9, for details. We need to make a adjustment for
$v_{2k}$ because of the singularity of $\log R$ at $R=0$. Also, we need to make sure that $v_{2k}$ has order 3 vanishing at $R= 0$. Thus we
define $v_{2k}$ as
\begin{align*}
&v_{2k}:=\\&
\frac{1}{(t\lambda)^{2k-1}}\sum_{j=0}^{2k-1}W^j_{2k}(a)\Big(\frac 12
\log (1+R^2)\Big)^j+\frac{1}{(t\lambda)^{2k}}\frac{R}{(1+R^2)^{\frac{1}{2}}}\sum_{j=0}^{2k}\widetilde{W}^j_{2k}(a)\Big(\frac 12
\log (1+R^2)\Big)^j
\end{align*}
By doing this we get a large error near $R=0$, but it is not very
significant since the purpose of the correction is to improve the
error for large $R$. Since $a=R/t\lambda$, it's easy to pull
out a $a^3$ factor from $W$'s and $a^2$ from $\widetilde{W}$'s to see
that we have \eqref{eq:v_2k}. 
\\

{\bf{Step 3}}: {\it{Show that the error $e_{2k}$ generated by $u_{2k} = u_{2k-1}+ v_{2k}$ satisfies \eqref{eq:e_2k}.}} Write 
$$t^2e_{2k}=t^2(e_{2k-1}-e^0_{2k-1})+t^2\big(e^0_{2k-1}-\tilde{\Box}(v_{2k})\big)+t^2N_{2k}(v_{2k})$$
where we recall \eqref{NLT}. We begin with the first term on the right, which has the form 
$$t^2(e_{2k-1}-e^0_{2k-1})\in \frac{1}{(t\lambda)^{2k}}[IS^1(R^{-1}(\log R)^{2k},\mathcal{Q}'_{k-1})+\langle b_1, b_2\rangle IS^1(R(\log R)^{2k-1},\mathcal{Q}'_{k-1})]$$
and we conclude by observing that 
$$IS^1(R^{-1}(\log R)^{2k},\mathcal{Q}'_{k-1})\subset IS^1(R^{-1}(\log R)^{2k},\mathcal{Q}_{k-1})+\langle b_1, b_2\rangle IS^1(R(\log R)^{2k-1},\mathcal{Q}'_{k-1})$$.
The reason for this is that for $w\in IS^1(R^{-1}(\log R)^{2k}, \mathcal{Q}'_{k-1})$ we  can write
$$w=(1-a^2)w+\frac {1}{(t\lambda)^2}R^2w.$$
For the second term in the definition of $t^2e_{2k}$, we have that by the same computation as in {\it{(2)}} of the preceding subsection 
\[
t^2\big(e^0_{2k-1}-\tilde{\Box}(v_{2k})\big)\in \frac{1}{(t\lambda)^{2k}}\big[IS^1\big(R^{-1}(\log R)^{2k}, \mathcal{Q}_{k}\big) + \langle b_1, b_2\rangle[IS^1\big(R(\log R)^{2k-1}, \mathcal{Q}'_{k}\big)
\]
Finally, for the contribution of $t^2N_{2k}(v_{2k})$,  we use as in \cite{KST0} that 
$$u_{2k-1}-u_0\in\frac{1}{(t\lambda)^2}IS^3(R\log R, \mathcal{Q}_k)$$
and, reasoning as in \cite{KST0}, we find 
\begin{align*}
&t^2\frac{1-\cos (2u_{2k-1})}{r^2}v_{2k}\\&\in\frac{1}{(t\lambda)^{2k}}\Big(IS^3(R^{-1}(\log R)^{2k-1},\mathcal{Q}_k)+\frac{\langle b_1, b_2\rangle}{(t\lambda)^2}IS^5(R(\log R)^{2k-1},\mathcal{Q}_k)\Big)\\
&t^2\frac{\sin (2u_{2k-1})}{2r^2}(1-\cos (2v_{2k}))\\&\in \frac{1}{(t\lambda)^{2k}}\Big(IS^1(R^{-1}(\log R)^{2k},\mathcal{Q}_k)+{\langle b_1, b_2\rangle}IS^3(R(\log R)^{2k-1},\mathcal{Q}_k)\Big)\\
&t^2\frac{\cos (2u_{2k-1})}{2r^2}(2v_{2k}-\sin (2v_{2k}))\in \frac{\langle b_1, b_2\rangle}{(t\lambda)^{2k}}IS^1(R(\log R)^{2k-1},\mathcal{Q}'_k)
\end{align*}
This shows that $e_{2k}$ has the desired form. 
\\

Iteration of {\bf{Step 1}} - {\bf{Step 3}} immediately furnishes the proof of Theorem~\ref{thm:approxsol} .

\section{Interlude: some spectral theory}

Here we quickly gather the spectral theory needed for the construction of the precise solution. This is a quick summary of results contained in \cite{KST0}. In the sequel, we shall often invoke the operator 
\[
\mathcal{L}: = -\frac{d^2}{dr^2} + \frac{3}{4r^2} - \frac{8}{(1+r^2)^2}
\]
acting on (a subspace of) $L^2(0, \infty)$. The actual domain is given by 
\[
\text{Dom}(\mathcal{L}) = \{f\in L^2(0, \infty):\,f,\,f'\in AC_{loc}\big((0,\infty)\big),\,\mathcal{L}f\in L^2\big((0,\infty)\big)
\]
The operator ${\mathcal{L}}$ is then self-adjoint with this domain. The spectrum $\text{spec}(\mathcal{L}) = [0,\infty)$ is purely absolutely continuous. We then have the following key results, identically stated and proved in \cite{KST0}: 

\begin{theorem}\label{thm:spectral1}(\cite{KST0})(a) For each $z\in \mathbf{C}$ there exists a fundamental system $\phi(r, z), \theta(r, z)$ for $\mathcal{L} - z$ which is analytic in $z$ for each $r>0$ and has the asymptotic behavior 
\[
\phi(r, z)\sim r^{\frac{3}{2}},\,\theta(r, z)\sim \frac{1}{2}r^{-\frac{1}{2}},\,\text{as}\,r\rightarrow 0.
\]
In particular, their Wronskian $W(\phi(\cdot, z), \theta(\cdot, z)) = 1$ for all $z\in \mathbf{C}$. Here $\phi(\cdot, z)$ is the Weyl-Titchmarsh solution of $\mathcal{L} - z$ at $r = 0$. The functions $\phi(\cdot, z)$, $\theta(\cdot, z)$ are real valued for $z\in \R$. 
\\
(b) For each $z\in {\mathbf{C}}$, $\Im z >0$, let $\psi^{+}(r, z)$ denote the Weyl-Titchmarsh solution of $\mathcal{L} - z$ at $r = +\infty$, normalized such that 
\[
\psi^+(r, z) \sim z^{-\frac{1}{4}}e^{iz^{\frac{1}{2}}r}\,\text{as}\,r\rightarrow +\infty, \Im z^{\frac{1}{2}}>0.
\]
If $\xi>0$, then the limit $\psi^+(r, \xi+ i0)$ exists point-wise for all $r>0$ and we denote it by $\psi^+(r, \xi)$. Moreover, define $\psi^-(\cdot, \xi) = \overline{\psi^+(\cdot, \xi)}$. Then $\psi^+(r, \xi), \psi^-(r, \xi)$ form a fundamental system of ${\mathcal{L}} - \xi$ with asymptotic behavior $\psi^{\pm}(r, \xi)\sim \xi^{-\frac{1}{4}}e^{\pm i\xi^{\frac{1}{2}}r}$ as $r\rightarrow\infty$. 
\\
(c) The spectral measure of $\mathcal{L}$ is absolutely continuous and its density is given by 
\[
\rho(\xi) = \frac{1}{\pi}\Im m(\xi+i0)\chi_{\xi>0}
\]
with the generalized Weyl-Titchmarsh function 
\[
m(\xi) = \frac{W(\theta(\cdot, \xi), \psi^+(\cdot, \xi))}{W(\psi^+(\cdot, \xi), \phi(\cdot, \xi))}.
\]
(d) The distorted Fourier transform defined as 
\[
\mathcal{F}: f\rightarrow \hat{f}(\xi): = \lim_{b\rightarrow\infty}\int_0^b \phi(r, \xi)f(r)\,dr
\]
is a unitary operator from $L^2(\R_{+})$ to $L^2(\R_{+}, \rho)$, and its inverse is given by 
\[
\mathcal{F}^{-1}: \hat{f}\rightarrow f(r) = \lim_{\mu\rightarrow\infty}\int_0^{\mu}\phi(r, \xi)\hat{f}(\xi)\rho(\xi)\,d\xi
\]
Here $\lim$ refers to the $L^2(\R_{+}, \rho)$, respectively the $L^2(\R_{+})$ limit. 
\end{theorem}

The next two propositions detail the precise analytic structure of the functions $\phi(r, z)$, $\psi^{\pm}(r, z)$. This will be pivotal for our arguments below. 

\begin{theorem}\label{thm:spectral2}(\cite{KST0}) For any $z\in \mathbf{C}$, the fundamental system $\phi(r, z)$, $\theta(r, z)$ from the preceding theorem admits absolutely convergent asymptotic expansions 
\[
\phi(r, z) = \phi_0(r) + r^{-\frac{1}{2}}\sum_{j=1}^\infty (r^2 z)^j\phi_j(r^2)
\]
\[
\theta(r, z) = r^{-\frac{1}{2}}\frac{1}{2}\big(1-r^2 - \sum_{j=1}^\infty(r^2z)^j\theta_j(r^2)\big) - (2+\frac{4}{z})\phi(r, z)\log r
\]
where the functions $\phi_j, \theta_j$ are holomorphic in $U = \{\Re u >-\frac{1}{2}\}$ and satisfy the bounds 
\[
|\phi_j(u)|\leq \frac{3 C^j}{(j-1)!}\log (1+|u|),\,|\phi_1(u)|>\frac{1}{2}\log u\,\,\,\text{if}\,\,\,u\gg 1
\]
\[
|\theta_1(u)|\leq C|u|,\,|\theta_j(u)|\leq \frac{C^j}{(j-1)!}(1+|u|),\,u\in U.
\]
Furthermore, 
\[
\phi_1(u) = -\frac{1}{4}\log u + \frac{1}{2} + O(u^{-1}\log^2 u)\,\,\,\text{as}\,\,\,u\rightarrow\infty,
\]
as well as 
\[
\phi_1(u) = -\frac{u}{8} + \frac{u^2}{12} + O(u^3)\,\,\,\text{as}\,\,\,u\rightarrow 0.
\]
\end{theorem}

As for the functions $\psi^{\pm}(r, z)$, they admit Hankel expansions at infinity, as follows: 

\begin{theorem}\label{thm:spectral3}(\cite{KST0}) For any $\xi>0$, the solution $\psi^+(\cdot, \xi)$ from Theorem~\ref{thm:spectral1} is of the form 
\[
\psi^+(r, \xi) = \xi^{-\frac{1}{4}}e^{ir\xi^{\frac{1}{2}}}\sigma(r\xi^{\frac{1}{2}}, r),\,r\xi^{\frac{1}{2}}\gtrsim 1,
\]
where $\sigma$ admits the asymptotic series approximation 
\[
\sigma(q, r) \sim \sum_{j=0}^\infty q^{-j}\psi_j^+(r),\,\psi_0^+ = 1,\,\psi_1^+ = \frac{3i}{8}+O(\frac{1}{1+r^2})
\]
with zero order symbols $\psi_j^+(r)$ that are analytic at infinity, 
\[
\sup_{r>0}|(r\partial_r)^k\psi_j^+(r)|<\infty
\]
in the sense that for all large integers $j_0$, and all indices $\alpha, \beta$, we have 
\[
\sup_{r>0}\big|(r\partial_r)^{\alpha}(q\partial_q)^{\beta}\big[\sigma(q, r) - \sum_{j=0}^{j_0}q^{-j}\psi_j^+(r)\big]\big|\leq c_{\alpha, \beta, j_0}q^{-j_0 - 1}
\]
for all $q>1$. 
\end{theorem}
Finally, the structure of the spectral measure is given by the following 
\begin{theorem}\label{thm:spectral4}(\cite{KST0}) 
(a) We have 
\[
\phi(r, \xi) = a(\xi)\psi^+(r, \xi) + \overline{a(\xi)\psi^+(r, \xi)},
\]
where $a$ is smooth, always nonzero, and has size 
\[
|a(\xi)|\sim -\xi^{\frac{1}{2}}\log\xi,\,\xi\ll 1,\,|a(\xi)|\sim \xi^{-\frac{1}{2}},\,\xi\gtrsim 1
\]
Moreover, it satisfies the symbol bounds 
\[
|(\xi\partial\xi)^k a(\xi)|\leq c_k |a(\xi)|,\,\forall \xi>0.
\]
(b) The spectral measure $\rho(\xi)d\xi$ has density 
\[
\rho(\xi) = \frac{1}{\pi}|a(\xi)|^{-2}
\]
and therefore satisfies 
\[
\rho(\xi)\sim \frac{1}{\xi\log^2\xi},\,\xi\ll 1,\,\rho(\xi)\sim \xi,\,\xi\gtrsim 1.
\]

\end{theorem}

\section{Construction of the precise solution}

Our point of departure here is the assumption that an approximate solution $u_{2k-1}$, $k\gg 1$, has been constructed, with a corresponding error term $e_{2k-1}$ which decays rapidly in the renormalized time $\tau: = \nu^{-1} t^{-\nu}$. Note that with respect to this time, we get 
\[
\lambda(\tau) : = \lambda(t(\tau)) = (\nu\tau)^{\frac{1+\nu}{\nu}}
\]
We also have the re-scaled variable $R = \lambda(\tau)r$. 
We shall assume that 
\[
|e_{2k-1}(t, r)|\lesssim \tau^{-N},\,r\leq t
\]
for some sufficiently large $N$, which is possible if we choose $k$ large enough. We shall also assume the fine structure of $e_{2k-1}$ as in section~\ref{sec:approxsol}, and more specifically as in \eqref{eq:e_2k-1}. 
We try to complete the approximate solution $u_{2k-1}$ to an exact solution $u = u_{2k-1} + \eps$, where $\eps$ solves the following equation, see (3.2) on p. 16 in \cite{KST0}, and \eqref{NLT} in section~\ref{sec:approxsol} for the definition of $N_{2k-1}$:
\begin{equation}\label{eq:tilepsmain}
\big(-(\partial_{\tau} + \frac{\lambda_{\tau}}{\lambda}R\partial_R)^2 + \frac{1}{4}(\frac{\lambda_{\tau}}{\lambda})^2 + \frac{1}{2}\partial_{\tau}(\frac{\lambda_{\tau}}{\lambda})\big)\tilde{\eps} - \mathcal{L}\tilde{\eps} = \lambda^{-2}R^{\frac{1}{2}}\big(N_{2k-1}(R^{-\frac{1}{2}}\tilde{\eps}) + e_{2k-1}\big)
\end{equation}
Our strategy is to formulate this equation in terms of the Fourier coefficients of $\tilde{\eps}$ with respect to the Fourier basis associated with $ \mathcal{L}$, the latter as in the preceding section, given by 
\[
\mathcal{L} =  -\partial_R^2 + \frac{3}{4R^2} - \frac{8}{(1+R^2)^2}
\]

Introduce the operator 
\[
\mathcal{K} = -\big(\frac{3}{2} + \frac{\eta\rho'(\eta)}{\rho(\eta)}\big)\delta_0(\xi - \eta) + \mathcal{K}_0,
\]
see (5.3) on p. 25 in \cite{KST0}. This operator is needed to describe the transition from \eqref{eq:tilepsmain} to the equivalent formulation in terms of the Fourier coefficients: 
\begin{equation}\label{eq:Fourier1}\begin{split}
\mathcal{F}\big(\partial_{\tau} + \frac{\lambda_{\tau}}{\lambda}R\partial_R\big)^2& = \big(\partial_{\tau} + \frac{\lambda_{\tau}}{\lambda}(-2\xi\partial_{\xi} + \mathcal{K})\big)^2\mathcal{F}\\
& =  \big(\partial_{\tau} - \frac{\lambda_{\tau}}{\lambda}2\xi\partial_{\xi}\big)^2\mathcal{F} + \big(\partial_{\tau} - \frac{\lambda_{\tau}}{\lambda}2\xi\partial_{\xi}\big)\frac{\lambda_{\tau}}{\lambda}\mathcal{K}\mathcal{F}\\
&+\frac{\lambda_{\tau}}{\lambda}\mathcal{K}\big(\partial_{\tau} - \frac{\lambda_{\tau}}{\lambda}2\xi\partial_{\xi}\big)\mathcal{F} + (\frac{\lambda_{\tau}}{\lambda})^2\mathcal{K}^2\mathcal{F}\\
&= \big(\partial_{\tau} - \frac{\lambda_{\tau}}{\lambda}2\xi\partial_{\xi}\big)^2\mathcal{F} + 2\frac{\lambda_{\tau}}{\lambda}\mathcal{K}\big(\partial_{\tau} - \frac{\lambda_{\tau}}{\lambda}2\xi\partial_{\xi}\big)\mathcal{F}\\
&+\partial_{\tau}(\frac{\lambda_{\tau}}{\lambda})\mathcal{K}\mathcal{F} - 2(\frac{\lambda_{\tau}}{\lambda})^2[\xi\partial_{\xi}, \mathcal{K}]\mathcal{F} + (\frac{\lambda_{\tau}}{\lambda})^2\mathcal{K}^2\mathcal{F}\\
\end{split}\end{equation}
To proceed further, we have to precisely understand the structure of the 'transference operator' $\mathcal{K}$. Make the 
\begin{defi} We call an operator $\tilde{\mathcal{K}}$ to be 'smoothing', provided it enjoys the mapping property 
\[
\tilde{\mathcal{K}}: L^{2, \alpha}_\rho\longrightarrow  L^{2, \alpha+\frac{1}{2}}_\rho\,\,\,\forall \alpha
\]
\end{defi}
For the definition of the weighted $L^2$-space $ L^{2, \alpha}_\rho$ see (5.15), p. 29, in \cite{KST0}. Specifically, we have 
\[
\|u\|_{L^{2, \alpha}_\rho}: = \big(\int_0^\infty|u(\xi)|^2\langle\xi\rangle^{2\alpha}\rho(\xi)\,d\xi\big)^{\frac{1}{2}}
\]
The preceding definition means that applying $\tilde{\mathcal{K}}$ we gain $1/2$ power of $\xi$ of decay as $\xi\rightarrow\infty$. 
\\
For future reference, we will also use the following notation: if 
\[
f(R) = \int_0^\infty \phi(R, \xi)x(\xi)\rho(\xi)\,d\xi,
\]
then we write 
\[
\|f\|_{H^{\alpha}_{\rho}}: = \big(\int_0^\infty x^2(\xi)\langle\xi\rangle^{2\alpha}\,\rho(\xi)\,d\xi\big)^{\frac{1}{2}}
\]

Now according to Proposition 5.2 in \cite{KST0}, both operators $\mathcal{K}_0, [\xi\partial_\xi, \mathcal{K}_0]$, are smoothing. This is not stated this way in the cited Proposition for the commutator term, but the same proof as for $\mathcal{K}_0$ yields the smoothing property for $ [\xi\partial_\xi, \mathcal{K}_0]$. 
Our strategy shall be to move terms involving a smoothing operator to the right hand side, and keep those terms without smoothing property on the left, building them implicitly into the parametrix. This procedure is different than that employed in \cite{KST0}, and mimics the strategy in \cite{KS}. 
\\
Write (see Theorem~\ref{thm:spectral1})
\[
\tileps(\tau, R) = \int_0^\infty \phi(R, \xi)x(\tau, \xi)\rho(\xi)\,d\xi
\]
whence $x(\tau, \xi) = (\mathcal{F}\tileps)(\tau, \xi)$. Then using $\mathcal{F}\big(\mathcal{L}\tilde{\eps}\big)(\tau, \xi) = \xi x(\tau, \xi)$, we get from \eqref{eq:tilepsmain} and \eqref{eq:Fourier1}(all functions are to be evaluated at $(\tau, \xi)$)
\begin{equation}\label{eq:Fourier2}\begin{split}
- \big(\partial_{\tau} - \frac{\lambda_{\tau}}{\lambda}2\xi\partial_{\xi}\big)^2 x -\xi x& =  2\frac{\lambda_{\tau}}{\lambda}\mathcal{K}\big(\partial_{\tau} - \frac{\lambda_{\tau}}{\lambda}2\xi\partial_{\xi}\big)x + (\frac{\lambda_{\tau}}{\lambda})^2\big[\mathcal{K}^2 - 2[\xi\partial_{\xi}, \mathcal{K}]\big]x\\
&-\big(\frac{1}{4}(\frac{\lambda_{\tau}}{\lambda})^2 + \frac{1}{2}\partial_{\tau}(\frac{\lambda_{\tau}}{\lambda})\big)x + \partial_{\tau}(\frac{\lambda_{\tau}}{\lambda})\mathcal{K}x\\
&+\lambda^{-2}\mathcal{F}\big[R^{\frac{1}{2}}\big(N_{2k-1}(R^{-\frac{1}{2}}\tilde{\eps}) + e_{2k-1}\big)\big]
\end{split}\end{equation}
Here we want to remove all linear terms that do not have the smoothing property from the right hand side, which forces us to modify the left hand side. Observe the identity 
\begin{align*}
\big(\partial_{\tau} - \frac{\lambda_{\tau}}{\lambda}[2\xi\partial_{\xi}+\frac{3}{2}+\frac{\rho'(\xi)\xi}{\rho(\xi)}]\big)^2& = \big(\partial_{\tau} - 2\frac{\lambda_{\tau}}{\lambda}\xi\partial_{\xi}\big)^2\\& -  \big(\partial_{\tau} - 2\frac{\lambda_{\tau}}{\lambda}\xi\partial_{\xi}\big)\frac{\lambda_{\tau}}{\lambda}[\frac{3}{2}+\frac{\rho'(\xi)\xi}{\rho(\xi)}]\\
&-\frac{\lambda_{\tau}}{\lambda}[\frac{3}{2}+\frac{\rho'(\xi)\xi}{\rho(\xi)}]\big(\partial_{\tau} - 2\frac{\lambda_{\tau}}{\lambda}\xi\partial_{\xi}\big)\\
&+[\frac{\lambda_{\tau}}{\lambda}]^2[\frac{3}{2}+\frac{\rho'(\xi)\xi}{\rho(\xi)}]^2
\end{align*}
It follows that we have the relation 
\begin{equation}\label{eq:Fourier3}\begin{split}
-\big(\partial_{\tau} - \frac{\lambda_{\tau}}{\lambda}[2\xi\partial_{\xi}+\frac{3}{2}+\frac{\rho'(\xi)\xi}{\rho(\xi)}]\big)^2x -\xi x& =  2\frac{\lambda_{\tau}}{\lambda}\mathcal{K}_0\big(\partial_{\tau} - \frac{\lambda_{\tau}}{\lambda}2\xi\partial_{\xi}\big)x\\& + (\frac{\lambda_{\tau}}{\lambda})^2\big[\mathcal{K}^2 - (\mathcal{K} - \mathcal{K}_0)^2 - 2[\xi\partial_{\xi}, \mathcal{K}_0]\big]x\\
&-\big(\frac{1}{4}(\frac{\lambda_{\tau}}{\lambda})^2 + \frac{1}{2}\partial_{\tau}(\frac{\lambda_{\tau}}{\lambda})\big)x + \partial_{\tau}(\frac{\lambda_{\tau}}{\lambda})\mathcal{K}_0x\\
&+\lambda^{-2}\mathcal{F}\big[R^{\frac{1}{2}}\big(N_{2k-1}(R^{-\frac{1}{2}}\tilde{\eps}) + e_{2k-1}\big)\big]
\end{split}\end{equation}
Here the linear expression 
\[
\big(\frac{1}{4}(\frac{\lambda_{\tau}}{\lambda})^2 + \frac{1}{2}\partial_{\tau}(\frac{\lambda_{\tau}}{\lambda})\big)x = \tau^{-2}\big(\frac{1}{4}(\frac{\nu+1}{\nu}\big)^2 - \frac{1}{2}\frac{\nu+1}{\nu}\big)x =: c\tau^{-2}x
\]
still doesn't have the smoothing property. However, $x$ has better decay than the source terms $e_{2k-1}$, and so we will gain smoothness once we apply the parametrix to this term. We shall therefore leave it on the right hand side. 
To simplify notation, introduce the operator
\[
\mathcal{D}_{\tau}: = \partial_{\tau} - \frac{\lambda_{\tau}}{\lambda}[2\xi\partial_{\xi}+\frac{3}{2}+\frac{\rho'(\xi)\xi}{\rho(\xi)}]
\]
Then we can finally formulate \eqref{eq:Fourier3} in the form 
\begin{equation}\label{eq:Fourier4}
\mathcal{D}_{\tau}^2x + \xi x = f,
\end{equation}
where we have 
\begin{equation}\label{eq:Fourier5}\begin{split}
-f = &2\frac{\lambda_{\tau}}{\lambda}\mathcal{K}_0\big(\partial_{\tau} - \frac{\lambda_{\tau}}{\lambda}2\xi\partial_{\xi}\big)x + (\frac{\lambda_{\tau}}{\lambda})^2\big[\mathcal{K}^2 - (\mathcal{K} - \mathcal{K}_0)^2 - 2[\xi\partial_{\xi}, \mathcal{K}_0]\big]x\\
&+\partial_{\tau}(\frac{\lambda_{\tau}}{\lambda})\mathcal{K}_0x+\lambda^{-2}\mathcal{F}\big[R^{\frac{1}{2}}\big(N_{2k-1}(R^{-\frac{1}{2}}\tilde{\eps}) + e_{2k-1}\big)\big] - c\tau^{-2}x
\end{split}\end{equation}
In order to solve \eqref{eq:Fourier4}, we first formally replace $\mathcal{D}_{\tau}$ by $\partial_{\tau}$ and reduce to the simpler model problem 
\begin{equation}\label{eq:Fourier6}
L_{\xi, \tau}x: = \partial_{\tau}^2x + \lambda^{-2}(\tau)\xi x = f,
\end{equation}
The extra factor $\lambda^{-2}(\tau)$ comes from a change of scale. 
Introduce the symbol $S(\tau, \sigma, \xi)$ as in \cite{KST0}, via the requirements
\[
\partial_{\tau}^2 S  + \lambda^{-2}(\tau)\xi S = 0,\,S(\tau, \tau, \xi) = 0,\,\partial_{\tau}S(\tau, \sigma, \xi)|_{\tau = \sigma} = -1.
\]
We commence by noting that it suffices to consider the case $\xi = 1$. In fact (see \cite{KST0}), 
\begin{lemma}\label{lem:scaling}
We have the scaling relation 
\[
S(\tau, \sigma, \xi) = \xi^{\frac{\nu}{2}}S(\tau\xi^{-\frac{\nu}{2}}, \sigma\xi^{-\frac{\nu}{2}}, 1)
\]
\end{lemma}
\begin{proof}
We verify that the expression on the right has the desired properties. This follows from 
\begin{align*}
&\partial_{\tau}^2\xi^{\frac{\nu}{2}}S(\tau\xi^{-\frac{\nu}{2}}, \sigma\xi^{-\frac{\nu}{2}}, 1) = \xi^{-\frac{\nu}{2}}(\partial_{\tau}^2S)(\tau\xi^{-\frac{\nu}{2}}, \sigma\xi^{-\frac{\nu}{2}}, 1)\\
&\tau^{-2-\frac{2}{\nu}}\xi \big(\xi^{\frac{\nu}{2}}S(\tau\xi^{-\frac{\nu}{2}}, \sigma\xi^{-\frac{\nu}{2}}, 1)\big) = \xi^{-\frac{\nu}{2}}(\tau\xi^{-\frac{\nu}{2}})^{-2-\frac{2}{\nu}}S(\tau\xi^{-\frac{\nu}{2}}, \sigma\xi^{-\frac{\nu}{2}}, 1),
\end{align*}
where we recall that $\lambda(\tau) \sim \tau^{\frac{1+\nu}{\nu}}$. 
\end{proof}
We now construct $S(\tau, \sigma, 1)$ via the following 
\begin{lemma} (a) Let $\nu\leq \frac{1}{2}$. Then there exist two analytic solutions $\phi_0, \phi_1$ of $L_{1, \tau}\phi_j = 0$, $j = 0,1$, which admit a Puiseux series type representation 
\[
\phi_j(\tau) = \sum_{k=0}^\infty c_{jk}\tau^{j-\frac{2k}{\nu}},\,c_{j0} = 1,\,j\in \{0,1\}.
\]
The series converges absolutely for any $\tau>0$. 
(b) There is a solution $\phi_2(\tau)$ for $L_{1,\tau}$ of the form 
\[
\phi_2(\tau) = \tau^{\frac{1}{2}+\frac{1}{2\nu}}e^{i\nu\tau^{-\frac{1}{\nu}}}[1+a(\tau^{\frac{1}{\nu}})]
\]
with $a(0) = 0$.
\end{lemma}
This is Lemma 7. 1 of \cite{KST0}. It implies the following key 
\begin{prop}(\cite{KST0}) There is a constant $C>0$ such that we have the bounds 
\[
|S(\tau, \sigma, \xi)|\lesssim \sigma (\frac{\sigma}{\tau})^C(1+\tau^{-\frac{2}{\nu}}\xi)^{-\frac{1}{2}},\,|\partial_{\tau}S(\tau, \sigma, \xi)|\lesssim  (\frac{\sigma}{\tau})^C
\]
\end{prop}

We can now write down the explicit solution of \eqref{eq:Fourier4}: 
\begin{lemma} The equation \eqref{eq:Fourier4} is formally solved by the following parametrix
\begin{equation}\label{eq:Sparametrix}
x(\tau, \xi) = \int_{\tau}^\infty \frac{\lambda^{\frac{3}{2}}(\tau)}{\lambda^{\frac{3}{2}}(\sigma)}\frac{\rho^{\frac{1}{2}}(\frac{\lambda^2(\tau)}{\lambda^2(\sigma)}\xi)}{\rho^{\frac{1}{2}}(\xi)}S(\tau, \sigma, \lambda^2(\tau)\xi)f(\sigma, \frac{\lambda^2(\tau)}{\lambda^2(\sigma)}\xi)\,d\sigma =:(Uf)(\tau, \xi)
\end{equation}
\end{lemma}
This is a simple direct verification, as in \cite{KS}. For us, we will need the mapping properties of this parametrix with respect to suitable Banach spaces. We have 
\begin{lemma}\label{lem:ParaBounds}Introducing the norm 
\[
\|f\|_{L_{\rho}^{2,\alpha;N}}: = \sup_{\tau>\tau_0}\tau^{N}\|f(\tau, \cdot)\|_{L^{2,\alpha}_{\rho}},
\]
we have 
\[
\|Uf\|_{L_{\rho}^{2,\alpha+\frac{1}{2};N-2}}\lesssim \|f\|_{L_{\rho}^{2,\alpha;N}}
\]
provided $N$ is sufficiently large. 
\end{lemma}
\begin{proof}  This is a consequence of the bounds in the preceding proposition. Observe that 
\[
\frac{\rho^{\frac{1}{2}}(\frac{\lambda^2(\tau)}{\lambda^2(\sigma)}\xi)}{\rho^{\frac{1}{2}}(\xi)}\lesssim \frac{\lambda(\sigma)}{\lambda(\tau)}
\]
It follows that 
\begin{align*}
&\big|\langle \xi\rangle^{\alpha+\frac{1}{2}} \frac{\lambda^{\frac{3}{2}}(\tau)}{\lambda^{\frac{3}{2}}(\sigma)}\frac{\rho^{\frac{1}{2}}(\frac{\lambda^2(\tau)}{\lambda^2(\sigma)}\xi)}{\rho^{\frac{1}{2}}(\xi)}S(\tau, \sigma, \lambda^2(\tau)\xi)f(\sigma, \frac{\lambda^2(\tau)}{\lambda^2(\sigma)}\xi)\big|\\
&\lesssim \langle \xi\rangle^{\alpha+\frac{1}{2}}\big(\frac{\lambda(\tau)}{\lambda(\sigma)}\big)^{\frac{1}{2}}\sigma (\frac{\sigma}{\tau})^C(1+\tau^2\xi)^{-\frac{1}{2}}|f(\sigma, \frac{\lambda^2(\tau)}{\lambda^2(\sigma)}\xi)|\\
&\lesssim \sigma (\frac{\sigma}{\tau})^C \big(\frac{\lambda(\sigma)}{\lambda(\tau)}\big)^{2\alpha - \frac{1}{2}}\big|\langle \frac{\lambda^2(\tau)}{\lambda^2(\sigma)}\xi\rangle^{\alpha}f(\sigma, \frac{\lambda^2(\tau)}{\lambda^2(\sigma)}\xi)\big|
\end{align*}
It follows that 
\begin{align*}
\big\|Uf\big\|_{L_{\rho}^{2,\alpha+\frac{1}{2};N-2}}&\lesssim\sup_{\sigma>\tau_0}\big\|f\big\|_{L^{2,\alpha; N}_{\rho}}\sup_{\tau>\tau_0}\tau^{N-2}\int_{\tau}^\infty\sigma (\frac{\sigma}{\tau})^C \big(\frac{\lambda(\sigma)}{\lambda(\tau)}\big)^{2\alpha + \frac{1}{2}}\sigma^{-N}\,d\sigma\\
&\lesssim \sup_{\sigma>\tau_0}\big\|f\big\|_{L^{2,\alpha; N}_{\rho}}
\end{align*}
provided $N>(2\alpha + \frac{1}{2})\frac{\nu+1}{\nu} + C + 2$. 

\end{proof}

The goal is now to formulate \eqref{eq:Fourier4}, \eqref{eq:Fourier5} as an integral equation and find a suitable fixed point, which will be the desired $x(\tau, \xi)$. 
The issue is that $x$ will only have very weak regularity a priori, in fact $x(\tau, \cdot)\in L_{\rho}^{2,\frac{1}{2}+\frac{\nu}{2}-}$ is optimal, see \cite{KST0}, and this does not suffice for good algebra estimates provided $\nu\leq \frac{1}{2}$. We thus have to find some better space to place $x$ into. The key for this will be the {\it{zeroth iterate}} for solving \eqref{eq:Fourier4}, \eqref{eq:Fourier5}.  Thus solve these via 
\begin{equation}\label{eq:fixedpoint}
x(\tau, \xi) = (Uf)(\tau, \xi)
\end{equation}
with $f$ as in \eqref{eq:Fourier5}. To find such a fixed point, we use the iterative scheme 
\[
x_j(\tau, \xi) = (Uf_{j-1})(\tau, \xi),\,j\geq 1
\]
where we put 
\begin{align*}
-f_j = &2\frac{\lambda_{\tau}}{\lambda}\mathcal{K}_0\big(\partial_{\tau} - \frac{\lambda_{\tau}}{\lambda}2\xi\partial_{\xi}\big)x_j + (\frac{\lambda_{\tau}}{\lambda})^2\big[\mathcal{K}^2 - (\mathcal{K} - \mathcal{K}_0)^2 - 2[\xi\partial_{\xi}, \mathcal{K}_0]\big]x_j\\
&+\partial_{\tau}(\frac{\lambda_{\tau}}{\lambda})\mathcal{K}_0x_j+\lambda^{-2}\mathcal{F}\big[R^{\frac{1}{2}}\big(N_{2k-1}(R^{-\frac{1}{2}}\tilde{\eps}_j) + e_{2k-1}\big)\big] - c\tau^{-2}x_j
\end{align*}
and of course we put 
\[
\tileps_j(\tau, R) = \int_0^\infty \phi(R, \xi)x_j(\tau, \xi)\rho(\xi)\,d\xi
\]
The zeroth iterate in turn is defined via 
\[
x_0(\tau, \xi) = (U\lambda^{-2}\mathcal{F}\big[R^{\frac{1}{2}}\big(e_{2k-1}\big)\big])(\tau, \xi);
\]
here we may also replace $e_{2k-1}$ by a function which co-incides with it in the backward light cone $r\leq t$, in light of Huyghen's principle. This shall be handy below. 
\subsection{The zeroth iterate}
We claim in effect the following: 
\begin{prop}\label{prop:key} There exists $\tilde{e}_{2k-1}\in H_{RdR}^{\frac{\nu}{2}-}$ such that $\tilde{e}_{2k-1}|_{r\leq t} = e_{2k-1}$, and such that if we put 
\[
x_0(\tau, \xi) = (U\lambda^{-2}\mathcal{F}\big[R^{\frac{1}{2}}\big(\tilde{e}_{2k-1}\big)\big])(\tau, \xi),
\]
we can write 
\[
x_0 = x_0^{(1)} + x_0^{(2)},
\]
where we have 
\[
x_0^{(1)}\in \tau^{-N}L_\rho^{2, \frac{1}{2}+\frac{\nu}{2}-}  
\]
as well as 
\[
\chi_{R<1}\tilde{\eps}_0^{(1)}(\tau, R) = \chi_{R<1}\int_0^\infty \phi(R, \xi)x_0^{(1)}(\tau, \xi)\rho(\xi)\,d\xi\in \tau^{-N}R^{\frac{3}{2}}L^{\infty}
\]
while also 
\[
x_0^{(2)}\in \tau^{-N}L_\rho^{2, 1+\frac{\nu}{2}-}  
\]
\end{prop}
\begin{remark}Note that for $R\geq 1$, we actually have the bound 
\[
\big|\tilde{\eps}_0^{(1)}(\tau, R)\big|\lesssim \tau^{-N}
\]
since $\tilde{\eps}_0^{(1)}(\tau, \cdot)\in H^{1+\nu}_{dR}$. 

\end{remark}
\begin{proof} Recall from structure of the error $e_{2k-1}$ of the approximate solution $u_{2k-1}$ that $e_{2k-1}$ can be written as a sum of terms involving the singular expressions 
\[
\tau^{-N-2}(1-a)^{i\nu-\frac{1}{2}}(\log (1-a))^j,\,j\leq j(i),\,i\geq 1,
\]
multiplied by smooth (in $t, r$) functions. In fact, up to an error of class $H_{\R^2}^{2+\nu-}$ (namely when $(2i-1)\nu>2+\nu$), there are only finitely many such expressions. 
We now define $\tilde{e}_{2k-1}$ by replacing each of the above expressions by their truncation 
\[
\tau^{-N-2}(1-a)^{i\nu-\frac{1}{2}}(\log (1-a))^j|_{r\leq t};
\]
and the rest of $e_{2k-1}$ is smoothly truncated to a dilate of the light cone $r\leq t$. Thus, specifically, we shall write
\[
\tilde{e}_{2k-1} = \tilde{e}^{(1)}_{2k-1} + \tilde{e}^{(2)}_{2k-1},
\]
where we may assume 
\[
\tilde{e}^{(2)}_{2k-1}\in \tau^{-N-2}H^{2+\nu-}_{\R^2}
\]
while $ \tilde{e}^{(1)}_{2k-1}$ is a sum of singular terms of the above form with smooth bounded coefficients. 
Since  $\mathcal{F}\circ T^{-1}(H^{\alpha}_{\R^2}) = L^{2,\frac{\alpha}{2}}_{\rho}$, where $T$ is the map 
\[
u(R)\rightarrow e^{i\theta}R^{-\frac{1}{2}}u(R),
\]
we see from Lemma~\ref{lem:ParaBounds} that we have the bound 
\[
\big\|U\lambda^{-2}\mathcal{F}\big[R^{\frac{1}{2}}\big(\tilde{e}^{(2)}_{2k-1}\big)\big]\big\|_{L^{2,1+\frac{\nu}{2}-; N}_{\rho}}\lesssim 1
\]
and so we can include $U\lambda^{-2}\mathcal{F}\big[R^{\frac{1}{2}}\big(\tilde{e}^{(2)}_{2k-1}\big)\big]$ into $x_0^{(2)}$.  
\\
Next, consider the more difficult contribution of $ \tilde{e}^{(1)}_{2k-1}$, where a more detailed analysis of the parametrix becomes necessary. 
For general $f$, consider the decomposition 
\begin{align*}
&\int_{\tau}^\infty \frac{\lambda^{\frac{3}{2}}(\tau)}{\lambda^{\frac{3}{2}}(\sigma)}\frac{\rho^{\frac{1}{2}}(\frac{\lambda^2(\tau)}{\lambda^2(\sigma)}\xi)}{\rho^{\frac{1}{2}}(\xi)}S(\tau, \sigma, \lambda^2(\tau)\xi)f(\sigma, \frac{\lambda^2(\tau)}{\lambda^2(\sigma)}\xi)\,d\sigma\\
&=\int_{\max\{\tau, C(\lambda^2(\tau)\xi)^{\frac{\nu}{2}}\}}^\infty \frac{\lambda^{\frac{3}{2}}(\tau)}{\lambda^{\frac{3}{2}}(\sigma)}\frac{\rho^{\frac{1}{2}}(\frac{\lambda^2(\tau)}{\lambda^2(\sigma)}\xi)}{\rho^{\frac{1}{2}}(\xi)}S(\tau, \sigma, \lambda^2(\tau)\xi)f(\sigma, \frac{\lambda^2(\tau)}{\lambda^2(\sigma)}\xi)\,d\sigma\\
&+\int_{\tau}^{\max\{\tau, C(\lambda^2(\tau)\xi)^{\frac{\nu}{2}}\}}\frac{\lambda^{\frac{3}{2}}(\tau)}{\lambda^{\frac{3}{2}}(\sigma)}\frac{\rho^{\frac{1}{2}}(\frac{\lambda^2(\tau)}{\lambda^2(\sigma)}\xi)}{\rho^{\frac{1}{2}}(\xi)}S(\tau, \sigma, \lambda^2(\tau)\xi)f(\sigma, \frac{\lambda^2(\tau)}{\lambda^2(\sigma)}\xi)\,d\sigma\\
&=: (Uf)^{(1)} + (Uf)^{(2)}
\end{align*}
for some large constant $C$. In the first integral, we have 
\[
\sigma \geq C(\lambda^2(\tau)\xi)^{\frac{\nu}{2}},
\]
whence we obtain 
\[
\xi\leq (\frac{\sigma}{C})^{\frac{2}{\nu}}\lambda^{-2}(\tau)
\]
It follows that 
\[
\big\|(Uf)^{(1)}\big\|_{L^{2, 1+\frac{\nu}{2}-; N}_{\rho}}\lesssim \big\|f\big\|_{L^{2, \frac{\nu}{2}-; N+\frac{2}{\nu}+C_1}_{\rho}}
\]
and so we have gained smoothness for this terms at the expense of temporal decay. 
It thus remains to consider the term $(Uf)^{(2)}$, which in fact requires most of the work. On account of Lemma~\ref{lem:scaling}, we have 
\[
S(\tau, \sigma, \lambda^2(\tau)\xi) = (\lambda^2(\tau)\xi)^{\frac{\nu}{2}}S(\tau(\lambda^2(\tau)\xi)^{-\frac{\nu}{2}}, \sigma(\lambda^2(\tau)\xi)^{-\frac{\nu}{2}}, 1)
\]
Then from the proof of Lemma 7. 1 in \cite{KST0}, we can write 
\begin{align*}
S(\tau(\lambda^2(\tau)\xi)^{-\frac{\nu}{2}}, \sigma(\lambda^2(\tau)\xi)^{-\frac{\nu}{2}}, 1) = \Im\big(\phi_2(\tau(\lambda^2(\tau)\xi)^{-\frac{\nu}{2}})\overline{\phi_2}( \sigma(\lambda^2(\tau)\xi)^{-\frac{\nu}{2}})\big)
\end{align*}
and so using the factorization\footnote{Here we correct a typo in \cite{KST0}: we replace a factor $\nu$ by the correct $\nu^{-\frac{1}{\nu}}$} $\phi_2(\tau) = \tau^{\frac{1}{2}+\frac{1}{2\nu}}e^{i\nu^{-\frac{1}{\nu}}\tau^{-\frac{1}{\nu}}}[1+a(\tau^{\frac{1}{\nu}})]$ as in Lemma 7.1 in \cite{KST0}, we obtain 
\begin{equation}\label{eq:Sasympto}\begin{split}
&(\lambda^2(\tau)\xi)^{\frac{\nu}{2}}S(\tau(\lambda^2(\tau)\xi)^{-\frac{\nu}{2}}, \sigma(\lambda^2(\tau)\xi)^{-\frac{\nu}{2}}, 1)\\& = \tau^{\frac{1}{2}+\frac{1}{2\nu}}\sigma^{\frac{1}{2}+\frac{1}{2\nu}}(\lambda^2(\tau)\xi)^{-\frac{1}{2}}\sin\big(\nu\xi^{\frac{1}{2}}\tau\big(1-(\frac{\tau}{\sigma})^{\frac{1}{\nu}}\big)\big)\big(1+a(\tau(\lambda^2(\tau)\xi)^{-\frac{\nu}{2}})\big)\big(1+a(\sigma(\lambda^2(\tau)\xi)^{-\frac{\nu}{2}})\big)\\
&=\big(\frac{\sigma}{\tau}\big)^{\frac{1}{2}+\frac{1}{2\nu}}\xi^{-\frac{1}{2}}\sin\big(\nu\xi^{\frac{1}{2}}\tau\big(1-(\frac{\tau}{\sigma})^{\frac{1}{\nu}}\big)\big)\big(1+a(\tau(\lambda^2(\tau)\xi)^{-\frac{\nu}{2}})\big)\big(1+a(\sigma(\lambda^2(\tau)\xi)^{-\frac{\nu}{2}})\big)\\
\end{split}\end{equation}
Here the function $a(\tau)$ is smooth with bounded derivatives. 
\\
Our task now consists in checking how the oscillations of this expression potentially cancel against the oscillations of $f(\sigma, \frac{\lambda^2(\tau)}{\lambda^2(\sigma)}\xi)$ in $(Uf)^{(2)}$. Recall that 
\[
f(\sigma, \frac{\lambda^2(\tau)}{\lambda^2(\sigma)}\xi) = \lambda^{-2}(\sigma)\mathcal{F}\big[R^{\frac{1}{2}}\big(\tilde{e}^{(1)}_{2k-1}(\sigma, \cdot)\big)\big](\frac{\lambda^2(\tau)}{\lambda^2(\sigma)}\xi))
\]
We need to analyze the large frequency asymptotics of this expression. We recall from Theorem~\ref{thm:spectral4} that 
\[
\phi(R, \xi) = a(\xi)\psi^+(R, \xi) + \overline{a(\xi)\psi^+(R, \xi)}
\]
where we have the large frequency asymptotics $|a(\xi)|\sim \xi^{-\frac{1}{2}}$, $\xi\gg 1$. The function $a(\xi)$ is smooth and in fact obeys symbol behavior, see Theorem~\ref{thm:spectral4}. 
Furthermore, the oscillatory function $\psi^+$ can be written in the form 
\[
\psi^+(R, \xi) = \xi^{-\frac{1}{4}}e^{iR\xi^{\frac{1}{2}}}\sigma(R\xi^{\frac{1}{2}}, R),\,R\xi^{\frac{1}{2}}\gtrsim 1,
\]
where we have a symbolic expansion, see Theorem~\ref{thm:spectral3},
\[
\sigma(q, r) = \sum_{j=0}^\infty q^{-j}\psi_j^+(r)
\]
and the functions $\psi_j^+$ are uniformly bounded and smooth with symbol behavior. 
We insert these asymptotics into the formula for the Fourier transform, using the singular source term 
\[
\lambda^{-2}R^{\frac{1}{2}}\tilde{e}^{(1)}_{2k-1} = \tau^{-N-2}\chi_{r\leq t}a^{\frac{1}{2}}(1-a)^{i\nu-\frac{1}{2}}(\log (1-a))^j,\,i\geq 1.
\]
In fact, we may replace all additional factors $R^{-k}(\log R)^l$ by $(\lambda(\sigma)\cdot\sigma)^{-k}(\log (\lambda(\sigma)\sigma))^l$, since the errors committed thereby gain one 
degree of smoothness, and are thus in $H^{1+\nu-}_{\R^2}$. By the same token, we can also include a smooth cutoff $\chi_{a\geq \frac{1}{2}}$.\\
We now find that (with $f(\sigma, \xi) = \mathcal{F}\big(\lambda^{-2}R^{\frac{1}{2}}\tilde{e}^{(1)}_{2k-1}(\sigma, \cdot)\big)(\xi)$ as well as $a = \frac{R}{\lambda(\sigma)\sigma}$)
\begin{align*}
f(\sigma, \xi) =& \sigma^{-N-2}\int_{0}^{\nu\sigma}a(\xi)\xi^{-\frac{1}{4}}e^{iR\xi^{\frac{1}{2}}}\sigma(R\xi^{\frac{1}{2}}, R)\chi_{a\geq \frac{1}{2}}a^{\frac{1}{2}}(1-a)^{i\nu-\frac{1}{2}}(\log (1-a))^j\,dR\\
&+ \sigma^{-N-2}\int_{0}^{\nu\sigma}\overline{a(\xi)}\xi^{-\frac{1}{4}}e^{-iR\xi^{\frac{1}{2}}}\overline{\sigma(R\xi^{\frac{1}{2}}, R)}\chi_{a\geq \frac{1}{2}}a^{\frac{1}{2}}(1-a)^{i\nu-\frac{1}{2}}(\log (1-a))^j\,dR\\
\end{align*}
Using the asymptotic expansion 
\[
\sigma(R\xi^{\frac{1}{2}}, R) = c_0 + O(\frac{1}{R\xi^{\frac{1}{2}}}),
\]
where the $O$-term enjoys symbol behavior, we get 
\begin{align*}
&\int_{0}^{\nu\sigma}a(\xi)\xi^{-\frac{1}{4}}e^{iR\xi^{\frac{1}{2}}}\sigma(R\xi^{\frac{1}{2}}, R)\chi_{a\geq \frac{1}{2}}a^{\frac{1}{2}}(1-a)^{i\nu-\frac{1}{2}}(\log (1-a))^j\,dR\\
&=c_0\int_{0}^{\nu\sigma}a(\xi)\xi^{-\frac{1}{4}}e^{iR\xi^{\frac{1}{2}}}\chi_{a\geq \frac{1}{2}}a^{\frac{1}{2}}(1-a)^{i\nu-\frac{1}{2}}(\log (1-a))^j\,dR\\
&+\int_{0}^{\nu\sigma}a(\xi)\xi^{-\frac{1}{4}}e^{iR\xi^{\frac{1}{2}}}O(R^{-1}\xi^{-\frac{1}{2}})\chi_{a\geq \frac{1}{2}}a^{\frac{1}{2}}(1-a)^{i\nu-\frac{1}{2}}(\log (1-a))^j\,dR\\
\end{align*}
To bound the second term, observe that 
\begin{align*}
&\int_{0}^{\nu\sigma}a(\xi)\xi^{-\frac{1}{4}}e^{iR\xi^{\frac{1}{2}}}O(R^{-1}\xi^{-\frac{1}{2}})\chi_{a\geq \frac{1}{2}}a^{\frac{1}{2}}(1-a)^{i\nu-\frac{1}{2}}(\log (1-a))^j\,dR\\
&=\int_{0}^{\nu\sigma}a(\xi)\xi^{-\frac{1}{4}}e^{iR\xi^{\frac{1}{2}}}O(R^{-\frac{1}{2}}\xi^{-\frac{1}{2}})(\nu\sigma)^{-\frac{1}{2}}\chi_{a\geq \frac{1}{2}}(1-a)^{i\nu-\frac{1}{2}}(\log (1-a))^j\,dR\\
&=O(\xi^{-\frac{7}{4}})
\end{align*}
after integration by parts with respect to $R$. 
In short, we have now shown that 
\begin{align*}
f(\sigma, \xi) = &c_0\int_{0}^{\nu\sigma}a(\xi)\xi^{-\frac{1}{4}}e^{iR\xi^{\frac{1}{2}}}\chi_{a\geq \frac{1}{2}}a^{\frac{1}{2}}(1-a)^{i\nu-\frac{1}{2}}(\log (1-a))^j\,dR\\
&+\overline{c_0\int_{0}^{\nu\sigma}a(\xi)\xi^{-\frac{1}{4}}e^{iR\xi^{\frac{1}{2}}}\chi_{a\geq \frac{1}{2}}a^{\frac{1}{2}}(1-a)^{i\nu-\frac{1}{2}}(\log (1-a))^j\,dR}\\
+O(\xi^{-\frac{7}{4}})
\end{align*}
We now analyze the integrals more closely. We introduce the variable $x = \nu\sigma - R$. Then we can write 
\begin{align*}
&\int_{0}^{\nu\sigma}a(\xi)\xi^{-\frac{1}{4}}e^{iR\xi^{\frac{1}{2}}}\chi_{a\geq \frac{1}{2}}a^{\frac{1}{2}}(1-a)^{i\nu-\frac{1}{2}}(\log (1-a))^j\,dR\\
&= e^{i\nu\sigma\xi^{\frac{1}{2}}}a(\xi)\xi^{-\frac{1}{4}}(\nu\sigma)^{\frac{1}{2}-i\nu}\int_0^{\infty}e^{ix\xi^{\frac{1}{2}}}\chi_{x\leq \frac{\nu\sigma}{2}}\big(1-\frac{x}{\nu\sigma}\big)^{\frac{1}{2}}x^{-\frac{1-2i\nu}{2}}\big(\log(\frac{x}{\nu\sigma})\big)^j\,dx
\end{align*}
Changing variables to $y = x\xi^{\frac{1}{2}}$ allows us to express this expression in the form 
\[
 e^{i\nu\sigma\xi^{\frac{1}{2}}}a(\xi)\xi^{-\frac{1}{2}-i\nu}F(\sigma, \xi),
\]
where we have 
\[
F(\sigma, \xi): = 
(\nu\sigma)^{\frac{1}{2}-i\nu}\int_0^\infty e^{iy}\chi_{y\leq  \frac{\nu\sigma\xi^{\frac{1}{2}}}{2}}\big(1-\frac{y}{\nu\sigma\xi^{\frac{1}{2}}}\big)^{\frac{1}{2}}y^{-\frac{1-2i\nu}{2}}\big(\log(\frac{y}{\nu\sigma\xi^{\frac{1}{2}}})\big)^j\,dy
\]
Observe that $F(\sigma, \xi)\in C^\infty$, and we have 
\[
\big|\partial_{\xi^{\frac{1}{2}}}^lF(\sigma, \xi)\big|\lesssim (\nu\sigma)^{\frac{1}{2}-i\nu}\xi^{-\frac{l}{2}},\,\big|\partial_{\sigma}F(\sigma, \xi)\big|\lesssim (\nu\sigma)^{\frac{1}{2}-i\nu}\sigma^{-1}.
\]
Here it is of course important that we have the restriction $y\leq  \frac{\nu\sigma\xi^{\frac{1}{2}}}{2}$. 
We thus now have the representation 
\begin{equation}\label{eq:fasympto}\begin{split}
f(\sigma, \xi) = &c_0\sigma^{-N} e^{i\nu\sigma\xi^{\frac{1}{2}}}a(\xi)\xi^{-\frac{1}{2}-i\nu}F(\sigma, \xi)\\
&+\overline{c_0\sigma^{-N} e^{i\nu\sigma\xi^{\frac{1}{2}}}a(\xi)\xi^{-\frac{1}{2}-i\nu}F(\sigma, \xi)}\\
&+\sigma^{-N}O(\xi^{-\frac{7}{4}})
\end{split}\end{equation}
where we keep in mind the restriction that $\xi>1$, as we only care about the large frequency case. We shall now use this in the context of $(Uf)^{(2)}$, see above, with 
\[
f = \lambda^{-2}(\sigma)\mathcal{F}\big[R^{\frac{1}{2}}\big(\tilde{e}^{(1)}_{2k-1}(\sigma, \cdot)\big](\xi)
\]
Begin by  writing 
\begin{align*}
(Uf)^{(2)}(\tau, \xi) =&\int_{\tau}^{\min\{C(\lambda^2(\tau)\xi)^{\frac{\nu}{2}}, \xi^{\frac{\nu}{2(1+\nu}}\tau\}}\frac{\lambda^{\frac{3}{2}}(\tau)}{\lambda^{\frac{3}{2}}(\sigma)}\frac{\rho^{\frac{1}{2}}(\frac{\lambda^2(\tau)}{\lambda^2(\sigma)}\xi)}{\rho^{\frac{1}{2}}(\xi)}S(\tau, \sigma, \lambda^2(\tau)\xi)f(\sigma, \frac{\lambda^2(\tau)}{\lambda^2(\sigma)}\xi)\,d\sigma\\
&+\int_{\min\{C(\lambda^2(\tau)\xi)^{\frac{\nu}{2}}, \xi^{\frac{\nu}{2(1+\nu}}\tau\}}^{\max\{\tau, C(\lambda^2(\tau)\xi)^{\frac{\nu}{2}}\}}\frac{\lambda^{\frac{3}{2}}(\tau)}{\lambda^{\frac{3}{2}}(\sigma)}\frac{\rho^{\frac{1}{2}}(\frac{\lambda^2(\tau)}{\lambda^2(\sigma)}\xi)}{\rho^{\frac{1}{2}}(\xi)}S(\tau, \sigma, \lambda^2(\tau)\xi)f(\sigma, \frac{\lambda^2(\tau)}{\lambda^2(\sigma)}\xi)\,d\sigma\\
&=: (Uf)^{(21)}(\tau, \xi) + (Uf)^{(22)}(\tau, \xi) 
\end{align*}
In the second integral, we have 
\[
\xi<\big(\frac{\sigma}{\tau}\big)^{\frac{2(1+\nu)}{\nu}}
\]
and so we get 
\[
\big\|(Uf)^{(22)}\big\|_{L_{\rho}^{2,1; N}}\lesssim \big\|f\big\|_{L_{\rho}^{2,\frac{\nu}{2}-; N-2}},
\]
provided $N$ is sufficiently large in relation to $\nu$. 
\\
We have now reduced things to $(Uf)^{(21)}(\tau, \xi)$, where we have $\frac{\lambda^2(\tau)}{\lambda^2(\sigma)}\xi>1$, and so we can use \eqref{eq:fasympto}. 
We shall combine this with the asymptotic relation \eqref{eq:Sasympto}. Just recording the integrand of the resulting expression and omitting constants, we find the schematic expression 
\begin{align*}
&\frac{\lambda^{\frac{3}{2}}(\tau)}{\lambda^{\frac{3}{2}}(\sigma)}\frac{\rho^{\frac{1}{2}}(\frac{\lambda^2(\tau)}{\lambda^2(\sigma)}\xi)}{\rho^{\frac{1}{2}}(\xi)}\big(\frac{\sigma}{\tau}\big)^{\frac{1}{2}+\frac{1}{2\nu}}\xi^{-\frac{1}{2}}\sin\big(\nu\xi^{\frac{1}{2}}\tau\big(1-(\frac{\tau}{\sigma})^{\frac{1}{\nu}}\big)\big)\prod_{\kappa = \tau, \sigma}\big(1+a(\kappa(\lambda^2(\tau)\xi)^{-\frac{\nu}{2}})\big)\\
&\hspace{3cm}\cdot \sigma^{-N} e^{\pm i\nu\sigma\frac{\lambda(\tau)}{\lambda(\sigma)}\xi^{\frac{1}{2}}}a(\frac{\lambda^2(\tau)}{\lambda^2(\sigma)}\xi)(\frac{\lambda^2(\tau)}{\lambda^2(\sigma)}\xi)^{-\frac{1}{2}- i\nu}F(\sigma, \frac{\lambda^2(\tau)}{\lambda^2(\sigma)}\xi)
\end{align*}
and so $(Uf)^{(21)}(\tau, \xi)$ is obtained by applying $\int_{\tau}^{\min\{C(\lambda^2(\tau)\xi)^{\frac{\nu}{2}}, \xi^{\frac{\nu}{2(1+\nu}}\tau\}}\,d\sigma$ to this integrand. 
Observe that the decay of this expression with respect to large $\xi$ is 
\[
O(\xi^{-\frac{3}{2} - i\nu}),
\]
but in order to obtain the desired $L^{2,1+\frac{\nu}{2}-; N}_{\rho}$-decay, we would need at least $\xi^{-2-\frac{\nu}{2}}$. The only way to eke out this extra decay in $\xi$ is to exploit the integration in $\sigma$, for which the product of the oscillatory factors 
\[
\sin\big(\nu\xi^{\frac{1}{2}}\tau\big(1-(\frac{\tau}{\sigma})^{\frac{1}{\nu}}\big)\big)\cdot e^{\pm i\nu\sigma\frac{\lambda(\tau)}{\lambda(\sigma)}\xi^{\frac{1}{2}}} = \frac{e^{i\big(\nu\xi^{\frac{1}{2}}\tau\big(1-(\frac{\tau}{\sigma})^{\frac{1}{\nu}}\big)\big)} - e^{-i\big(\nu\xi^{\frac{1}{2}}\tau\big(1-(\frac{\tau}{\sigma})^{\frac{1}{\nu}}\big)\big)}}{2}\cdot e^{\pm i\nu\sigma\frac{\lambda(\tau)}{\lambda(\sigma)}\xi^{\frac{1}{2}}}
\]
is key. The resulting phase functions (upon developing this product) are either of the form 
\[
e^{\pm i\big(\nu\xi^{\frac{1}{2}}\tau\big(1-2(\frac{\tau}{\sigma})^{\frac{1}{\nu}}\big)\big)},
\]
in which case we gain a factor $\xi^{-\frac{1}{2}}$ via integration by parts with respect to $\sigma$, or else of the form 
\[
e^{\pm i\nu\xi^{\frac{1}{2}}\tau},
\]
in which case the $\sigma$-oscillation has been destroyed. 
\\

It is this last case we now investigate more closely. We shall essentially put 
\[
x_0^{(1)} = (Uf)^{(21)}(\tau, \xi)
\]
Then the required inclusion $x_0^{(1)}\in L^{2, \frac{1}{2}+\frac{\nu}{2}-; N}_{\rho}$ is immediate, and so we now need to verify the sufficient vanishing of $\tileps_0^{(1)}(\tau, R)$ at $R = 0$. Thus consider 
\begin{align}
\tileps_0^{(1)}(\tau, R) = &\int_0^\infty \phi(R, \xi)x_0^{(1)}(\tau, \xi)\rho(\xi)\,d\xi\nonumber\\
& = \int_0^\infty \chi_{\xi<1}\phi(R, \xi)x_0^{(1)}(\tau, \xi)\rho(\xi)\,d\xi\label{eq:smallfreq}\\
&+ \int_0^\infty \chi_{1\leq \xi<R^{-2}}\phi(R, \xi)x_0^{(1)}(\tau, \xi)\rho(\xi)\,d\xi\label{eq:medfreq}\\
&+ \int_0^\infty \chi_{\xi\geq R^{-2}}\phi(R, \xi)x_0^{(1)}(\tau, \xi)\rho(\xi)\,d\xi\label{eq:largefreq}
\end{align}
We have included smooth cutoffs to dilates of the indicated regions. 
Here the first term \eqref{eq:smallfreq} clearly is in $L^{2,1;N}_{\rho}$ and hence negligible. It remains to control the other two terms, for which we use the asymptotic expansions of $\phi(R, \xi)$. For the last term, use 
\[
\phi(R, \xi) = a(\xi)\psi^+(R, \xi) + \overline{a(\xi)\psi^+(R, \xi)}
\]
with
\[
\psi^+(R, \xi) = \xi^{-\frac{1}{4}}e^{iR\xi^{\frac{1}{2}}}\sigma(R\xi^{\frac{1}{2}}, R),\,R\xi^{\frac{1}{2}}\gtrsim 1,
\]
as well as $|a(\xi)|\lesssim \xi^{-\frac{1}{2}}$. Then keeping in mind the structure of $x_0^{(1)} = (Uf)^{(21)}(\tau, \xi)$, we can write (schematically)
\begin{align}
& \int_0^\infty \chi_{\xi\geq R^{-2}}\phi(R, \xi)x_0^{(1)}(\tau, \xi)\rho(\xi)\,d\xi\nonumber\\
&=\int_0^\infty a(\xi)\xi^{-\frac{1}{4}}\chi_{\xi\geq R^{-2}}e^{i[R\xi^{\frac{1}{2}}\pm\nu\xi^{\frac{1}{2}}\tau]}\sigma(R\xi^{\frac{1}{2}}, R)\big(\int_\tau^{\kappa(\tau, \xi)}G_1(\sigma, \tau, \xi)\,d\sigma\big)\rho(\xi)\,d\xi\label{eq:G_1}\\
&+\int_0^\infty a(\xi)\xi^{-\frac{1}{4}}\chi_{\xi\geq R^{-2}}e^{iR\xi^{\frac{1}{2}}}\sigma(R\xi^{\frac{1}{2}}, R)\big(\int_\tau^{\kappa(\tau, \xi)}G_2(\sigma, \tau, \xi)\,d\sigma\big)\rho(\xi)\,d\xi\label{eq:G_2}
\end{align}
where we have used the notation 
\[
\kappa(\tau, \xi) = \min\{C(\lambda^2(\tau)\xi)^{\frac{\nu}{2}}, \xi^{\frac{\nu}{2(1+\nu}}\tau\}
\]
as well as 
\begin{align*}
G_1(\sigma, \tau, \xi) =& \frac{\lambda^{\frac{3}{2}}(\tau)}{\lambda^{\frac{3}{2}}(\sigma)}\frac{\rho^{\frac{1}{2}}(\frac{\lambda^2(\tau)}{\lambda^2(\sigma)}\xi)}{\rho^{\frac{1}{2}}(\xi)}\big(\frac{\sigma}{\tau}\big)^{\frac{1}{2}+\frac{1}{2\nu}}\xi^{-\frac{1}{2}}\prod_{\kappa = \tau, \sigma}\big(1+a(\kappa(\lambda^2(\tau)\xi)^{-\frac{\nu}{2}})\big)\\&\cdot \sigma^{-N} a(\frac{\lambda^2(\tau)}{\lambda^2(\sigma)}\xi)(\frac{\lambda^2(\tau)}{\lambda^2(\sigma)}\xi)^{-\frac{1}{2}-i\nu}F(\sigma, \frac{\lambda^2(\tau)}{\lambda^2(\sigma)}\xi)
\end{align*}
Further, for the oscillatory second integral, we have 
\begin{align*}
G_2(\sigma, \tau, \xi) =&\frac{\lambda^{\frac{3}{2}}(\tau)}{\lambda^{\frac{3}{2}}(\sigma)}\frac{\rho^{\frac{1}{2}}(\frac{\lambda^2(\tau)}{\lambda^2(\sigma)}\xi)}{\rho^{\frac{1}{2}}(\xi)}\big(\frac{\sigma}{\tau}\big)^{\frac{1}{2}+\frac{1}{2\nu}}\xi^{-\frac{1}{2}}e^{i\big(\pm\nu\xi^{\frac{1}{2}}\tau\big(1-2(\frac{\tau}{\sigma})^{\frac{1}{\nu}}\big)\big)}\prod_{\kappa = \tau, \sigma}\big(1+a(\kappa(\lambda^2(\tau)\xi)^{-\frac{\nu}{2}})\big)\\
&\hspace{3cm}\cdot \sigma^{-N} \xi^{\frac{1}{2}}a(\frac{\lambda^2(\tau)}{\lambda^2(\sigma)}\xi)(\frac{\lambda^2(\tau)}{\lambda^2(\sigma)}\xi)^{-\frac{1}{2}- i\nu}F(\sigma, \frac{\lambda^2(\tau)}{\lambda^2(\sigma)}\xi)
 \end{align*}
The idea now is that in the first integral \eqref{eq:G_1}, we can perform an integration by parts with respect to $\xi^{\frac{1}{2}}$, provided the phase $R\pm\nu\tau$ is large, which is certainly the case if we restrict to $R<\frac{\nu\tau}{2}$. More precisely, this becomes possible once we split the $\xi$-integral into two, where the limit $\kappa(\tau, \xi)$ is a smooth function of $\xi$. 
Observe that 
\[
\big|G_1(\sigma, \tau, \xi)\big|\lesssim \Lambda(\sigma, \tau)\xi^{-\frac{3}{2}}
\]
for a suitable $\Lambda(\sigma, \tau)$. Performing an integration by parts with respect to $\xi^{\frac{1}{2}}$ in \eqref{eq:G_1} and assuming $N$ to be large enough (in relation to $\nu$), as well as using the bound $\chi_{\xi\geq R^{-2}}\xi^{-\frac{3}{4}}\lesssim R^{\frac{3}{2}}$, we then find 
\begin{align*}
&\big|\chi_{R<\frac{\nu\tau}{2}}\int_0^\infty a(\xi)\xi^{-\frac{1}{4}}\chi_{\xi\geq R^{-2}}e^{i[R\xi^{\frac{1}{2}}\pm\nu\xi^{\frac{1}{2}}\tau]}\sigma(R\xi^{\frac{1}{2}}, R)\big(\int_\tau^{\kappa(\tau, \xi)}G_1(\sigma, \tau, \xi)\,d\sigma\big)\rho(\xi)\,d\xi\big|\\
&\lesssim \tau^{-(N-1)}R^{\frac{3}{2}}
\end{align*}
Next, consider the integral \eqref{eq:G_2}. Here we perform the integration by parts inside the $\sigma$-integral, due to the oscillatory phase 
\[
e^{i\big(\pm\nu\xi^{\frac{1}{2}}\tau\big(1-2(\frac{\tau}{\sigma})^{\frac{1}{\nu}}\big)\big)}
\]
Indeed, we have 
\begin{align*}
e^{i\big(\pm\nu\xi^{\frac{1}{2}}\tau\big(1-2(\frac{\tau}{\sigma})^{\frac{1}{\nu}}\big)\big)}
= \mp\big(\frac{\sigma}{\tau}\big)^{1+\nu^{-1}}(2i\xi^{\frac{1}{2}})^{-1}\partial_{\sigma}\big(e^{i\big(\pm\nu\xi^{\frac{1}{2}}\tau\big(1-2(\frac{\tau}{\sigma})^{\frac{1}{\nu}}\big)\big)}\big)
\end{align*}
and so we gain one inverse power $\xi^{-\frac{1}{2}}$ at the expense of a weight $\big(\frac{\sigma}{\tau}\big)^{1+\nu^{-1}}$, and this is enough to force absolute integrability with respect to $\xi$ since $\rho(\xi)\sim \xi$ for large $\xi$. 
It follows that 
\begin{align*}
&\big|\int_0^\infty a(\xi)\xi^{-\frac{1}{4}}\chi_{\xi\geq R^{-2}}e^{iR\xi^{\frac{1}{2}}}\sigma(R\xi^{\frac{1}{2}}, R)\big(\int_\tau^{\kappa(\tau, \xi)}G_2(\sigma, \tau, \xi)\,d\sigma\big)\rho(\xi)\,d\xi\big|\\
&\lesssim \tau^{-(N-1)}R^{\frac{3}{2}},
\end{align*}
even irrespective of the size of $R$. This concludes the estimate for the term \eqref{eq:largefreq}.
\\

It remains to deal with \eqref{eq:medfreq}, where we use the expansion 
\[
\phi(R, \xi) = \phi_0(R) + R^{-\frac{1}{2}}\sum_{j=1}^\infty (R^2\xi)^j\phi_j(R^2),
\]
where the functions $\phi_j$ are smooth with very good bounds: 
\[
|\phi_j(u)|\leq \frac{3C^j}{(j-1)!}\log(1+|u|),
\]
see Theorem~\ref{thm:spectral3}. Then as in \eqref{eq:G_1}, \eqref{eq:G_2}, we decompose 
\begin{align}
& \int_0^\infty \chi_{\xi<R^{-2}}\phi(R, \xi)x_0^{(1)}(\tau, \xi)\rho(\xi)\,d\xi\nonumber\\
&=\int_0^\infty \chi_{\xi<R^{-2}}\phi(R, \xi)e^{i\nu\xi^{\frac{1}{2}}\tau}\big(\int_\tau^{\kappa(\tau, \xi)}G_1(\sigma, \tau, \xi)\,d\sigma\big)\rho(\xi)\,d\xi\label{eq:H_1}\\
&+\int_0^\infty \chi_{\xi<R^{-2}}\phi(R, \xi)\big(\int_\tau^{\kappa(\tau, \xi)}G_2(\sigma, \tau, \xi)\,d\sigma\big)\rho(\xi)\,d\xi\label{eq:H_2}
\end{align}
In the first integral on the right, we perform integration by parts with respect to $\xi^{\frac{1}{2}}$, gaining a factor $\tau^{-1}$. If the derivative falls on the function $\phi(R, \xi)$, we obtain the differentiated series 
\[
\sum_{j=1}^\infty j(R^2\xi)^{j-1}R^{\frac{3}{2}}\xi^{\frac{1}{2}}\phi_j(R^2)
\]
which is bounded in absolute value by 
\[
\big|\sum_{j=1}^\infty j(R^2\xi)^{j-1}R^{\frac{3}{2}}\xi^{\frac{1}{2}}\phi_j(R^2)\big|\lesssim R^{\frac{3}{2}}\log(2+R)
\]
When the derivative falls on the inner integral, the bound is the same as before, and the last integral \eqref{eq:H_2} is also bounded just like \eqref{eq:G_2}. This concludes the proof of the proposition. 
\end{proof}

For later reference, we need somewhat more refined information, which however easily follows from the preceding proof. We mention 
\begin{cor}\label{cor:refine} Denote by $P_{\lambda}$ the frequency localizers 
\[
\mathcal{F}\big(P_{<\lambda} f\big)(\xi) = \chi_{<\lambda}(\xi)\big(\mathcal{F}f\big)(\xi)
\]
where $ \chi_{<\lambda}(\xi)$ is a smooth cutoff function localizing to $\xi\lesssim \lambda$, as in \cite{KST0}; here $\lambda$ is a dyadic number. Then we have 
\[
\chi_{R<1}P_{<\lambda}\tileps_0^{(1)}\in \tau^{-N}R^{\frac{3}{2}}L^\infty
\]
uniformly in $\lambda>1$. Furthermore, for any integer $l\geq 0$, we have 
\[
\nabla_R^l R^{-\frac{3}{2}}P_{<\lambda}\tileps_0^{(1)} = O(\tau^{-N})
\]
uniformly in $\lambda>1$. 
\end{cor}

\subsection{Analysis of the nonlinear source terms}

From \eqref{NLT}, we recall the following formula for the main source term:
\begin{align}
&\lambda^{-2}R^{\frac{1}{2}}N_{2k-1}(R^{-\frac{1}{2}}\tileps) = \nonumber\\
&\frac{4\sin(u_0 - u_{2k})\sin(u_0 + u_{2k})}{R^2}\tileps\label{eq:N_1}\\
&+\frac{\sin(2u_{2k})}{2R^{\frac{3}{2}}}\big(1-\cos(2R^{-\frac{1}{2}}\tileps)\big)\label{eq:N_2}\\
&+\frac{\cos(2u_{2k})}{2R^{\frac{3}{2}}}\big(2R^{-\frac{1}{2}}\tileps - \sin(2R^{-\frac{1}{2}}\tileps)\big)\label{eq:N_3}
\end{align}
 According to the preceding proposition, we have 
 \[
 x_0\in \tau^{-N}L^{2, \frac{1}{2}+\frac{\nu}{2}-}_{\rho}
 \]
 whence 
 \[
 \tileps_0(\tau, \cdot)\in \tau^{-N}H^{\frac{1}{2}+\frac{\nu}{2}-}_{\rho}
 \]
 This means that for the source terms, we need at least $H^{\frac{\nu}{2}-}_{\rho}$-regularity. In fact, we can do much better for the term \eqref{eq:N_1}. Recall that 
 \[
 u_{2k} = u_0 + \sum_{j=1}^{2k}v_{j}
 \]
 where we have 
 \[
 v_{2j-1}\in \frac{1}{(t\lambda)^{2j}}IS^3\big(R(\log R)^{2j-1}, \mathcal{Q}_{j-1}),\,v_{2j} \in  \frac{1}{(t\lambda)^{2j+2}}IS^3\big(R(\log R)^{2j-1}, \mathcal{Q}_{j})
\]
This implies in particular that 
\[
\frac{\sin(u_0 -u_{2k})}{R}\in (\lambda t)^{-2}IS(\log R, \mathcal{Q}),\,\frac{\sin(u_0 +u_{2k})}{R}\in IS(R^{-1}, \mathcal{Q})
\]
Then we recall Lemma 8.1 from \cite{KST0}: 

\begin{lemma}\label{lem:8.1}(\cite{KST0}) Assume $|\alpha|<\frac{\nu}{2}+\frac{3}{4}$, $f\in IS(1, \mathcal{Q})$. Then we have 
\[
\big\| gf\big\|_{H^{\alpha}_{\rho}}\lesssim \|f\|_{H^{\alpha}_{\rho}}
\]
\end{lemma}

Application of this lemma yields the bound 
\begin{equation}\label{eq:N_1bound}
\big\|\frac{4\sin(u_0 - u_{2k})\sin(u_0 + u_{2k})}{R^2}\tileps\big\|_{H^{\frac{1}{2}+\frac{\nu}{2}-}_{\rho}}\lesssim (\lambda t)^{-2}\big\|\tileps\big\|_{H^{\frac{1}{2}+\frac{\nu}{2}-}_{\rho}}
\end{equation}

To deal with the truly nonlinear source terms \eqref{eq:N_2} and \eqref{eq:N_3}, we need the following multilinear estimates: 

\begin{lemma}\label{lem:multilin1} 
Assume $f,g\in H^{\frac{1}{2}+\frac{\nu}{2}-}_{\rho}\cap R^{\frac{3}{2}}L^\infty$, $P_{<\lambda}f, P_{<\lambda}g\in \log\lambda R^{\frac{3}{2}}L^\infty$ uniformly in $\lambda>1$. If also $\chi_{R<1}\nabla^l\big(R^{-\frac{3}{2}}P_{<\lambda}f\big)\in L^\infty$ uniformly in $\lambda>1$, $l\geq 0$, then we have 
\[
R^{-\frac{3}{2}}fg\in H^{\frac{1}{2}+\frac{\nu}{2}-\delta-}_{\rho}\cap R^{\frac{3}{2}}L^\infty.
\]
for arbitrarily small $\delta\in (0,\frac{\nu}{100}]$ (with implicit constant depending on $\delta$), and we also have 
\[
R^{-\frac{3}{2}}P_{<\lambda}\big(R^{-\frac{3}{2}}fg\big)\in \log \lambda L^\infty,\,R^{-1}P_{<\lambda}\big(R^{-\frac{3}{2}}fg\big)\in L^\infty
\]
uniformly in $\lambda>1$. If $f\in H^{\frac{1}{2}+\frac{\nu}{2}-\delta-}_{\rho}\cap R^{\frac{3}{2}}L^\infty$, $P_{<\lambda}f\in \log\lambda R^{\frac{3}{2}}L^\infty$ uniformly in $\lambda$, but $g\in H_{\rho}^{1+\frac{\nu}{2}-2\delta-}$, $\delta\in (0,\frac{\nu}{100}]$, then 
\[
R^{-\frac{3}{2}}fg\in H^{\frac{1}{2}+\frac{\nu}{2}-}_{\rho}\cap R^{\frac{3}{2}}L^\infty,\,R^{-\frac{3}{2}}P_{<\lambda}\big(R^{-\frac{3}{2}}fg\big)\in \log \lambda R^{\frac{3}{2}}L^\infty,\,R^{-1}P_{<\lambda}\big(R^{-\frac{3}{2}}fg\big)\in L^\infty
\]
The same conclusion obtains if both $f, g\in H^{1+\frac{\nu}{2}-2\delta-}_{\rho}$. 
Further, if $f, g\in (H^{\frac{1}{2}+\frac{\nu}{2}-\delta-}_{\rho}\cap R^{\frac{3}{2}}L^\infty)$,  as well as 
\[
P_{<\lambda}f\in RL^\infty,\,P_{<\lambda}g\in  RL^\infty,\,\chi_{R<1}\nabla_R^l\big(R^{-1}P_{<\lambda}f\big)\in L^\infty,\,l\geq 0,
\]
uniformly in $\lambda>1$,  or else one of $f,g\in H_{\rho}^{1+\frac{\nu}{2}-2\delta-}$,
we get for $j = 0,1$
\[
R^{-j}fg\in H^{\frac{1}{2}+\frac{\nu}{2}-\delta-}_{\rho}\cap RL^\infty,\,P_{<\lambda}\big(R^{-j}fg\big)\in RL^\infty,
\]
the latter inclusion uniformly in $\lambda>1$. 
\end{lemma}
\begin{proof}
Throughout $\lambda_{1,2}, \sigma$ are dyadic numbers. We mimic the proof of Lemma 8.5 in \cite{KST0}. Write 
\[
R^{-\frac{3}{2}}fg = \sum_{\lambda_{1,2}}\sum_{\sigma<\max\{\lambda_{1,2}\}}P_{\sigma}\big(R^{-\frac{3}{2}}P_{\lambda_1}f P_{\lambda_2}g\big) + \sum_{\lambda_{1,2}}\sum_{\sigma\geq \max\{\lambda_{1,2}\}}P_{\sigma}\big(R^{-\frac{3}{2}}P_{\lambda_1}f P_{\lambda_2}g\big) 
\]
To bound the first term, write 
\begin{equation}\label{eq:smallsigma}\begin{split} 
&\sum_{\lambda_{1,2}}\sum_{\sigma<\max\{\lambda_{1,2}\}}P_{\sigma}\big(R^{-\frac{3}{2}}P_{\lambda_1}f P_{\lambda_2}g\big)\\
&= \sum_{\lambda_{1}<\lambda_2}\sum_{\sigma<\max\{\lambda_{1,2}\}}P_{\sigma}\big(R^{-\frac{3}{2}}P_{\lambda_1}f P_{\lambda_2}g\big)\\
&+ \sum_{\lambda_{1}\geq \lambda_2}\sum_{\sigma<\max\{\lambda_{1,2}\}}P_{\sigma}\big(R^{-\frac{3}{2}}P_{\lambda_1}f P_{\lambda_2}g\big)  
\end{split}\end{equation}
Then we get for the first term (after summing over $\lambda_1$ only)
\begin{align*}
\sigma^{\frac{1}{2}+\frac{\nu}{2}-\delta-}\big\|R^{-\frac{3}{2}}P_{<\lambda_2}f P_{\lambda_2}g\big\|_{L^2}&\leq \sigma^{\frac{1}{2}+\frac{\nu}{2}-\delta-}\big\|R^{-\frac{3}{2}}P_{<\lambda_2}f \big\|_{L^\infty}\big\|P_{\lambda_2}g\big\|_{L^{2}}\\
&\lesssim  \big(\frac{\sigma}{\lambda_2}\big)^{\frac{1}{2}+\frac{\nu}{2}-\delta-}\lambda_2^{-\delta}\big\|R^{-\frac{3}{2}}P_{<\lambda_2}f \big\|_{L^\infty}\big\|P_{\lambda_2}g\big\|_{H_{\rho}^{\frac{1}{2}+\frac{\nu}{2}-}}
\end{align*}
which is more than acceptable in the case $\sigma<\lambda_2$ (allowing for square summation over $\sigma, \lambda_2$), even taking into account the logarithmic loss from the factor $\big\|R^{-\frac{3}{2}}P_{<\lambda_2}f \big\|_{L^\infty}$ on the right, thus controlling the first term on the right of \eqref{eq:smallsigma} in case $g\in H^{\frac{1}{2}+\frac{\nu}{2}-}_{\rho}$. 
If on the other hand $g\in H^{1+\frac{\nu}{2}-2\delta-}_{\rho}$, we get 
\begin{align*}
\sigma^{\frac{1}{2}+\frac{\nu}{2}-}\big\|R^{-\frac{3}{2}}P_{<\lambda_2}f P_{\lambda_2}g\big\|_{L^2}&\leq \sigma^{\frac{1}{2}+\frac{\nu}{2}-}\big\|R^{-\frac{3}{2}}P_{<\lambda_2}f \big\|_{L^\infty}\big\|P_{\lambda_2}g\big\|_{L^{2}}\\
&\leq \big(\frac{\sigma}{\lambda_2}\big)^{\frac{1}{2}+\frac{\nu}{2}-}\lambda_2^{-\frac{1}{2}+2\delta}\big\|R^{-\frac{3}{2}}P_{<\lambda_2}f \big\|_{L^\infty}\big\|P_{\lambda_2}g\big\|_{H_{\rho}^{1+\frac{\nu}{2}-2\delta-}}\\
\end{align*}
Here, we can again square-sum over $\sigma, \lambda_2$. 
Next, for the case $\lambda_1\geq \lambda_2$ in \eqref{eq:smallsigma}, the argument is identical to the one above provided $g\in H^{\frac{1}{2}+\frac{\nu}{2}-}_{\rho}\cap R^{\frac{3}{2}}L^\infty$, $P_{<\lambda}f,g\in \log\lambda R^{\frac{3}{2}}L^\infty$ uniformly in $\lambda>1$. On the other hand, if $g\in H^{1+\frac{\nu}{2}-2\delta-}_{\rho}$, we have
\begin{align*}
\sigma^{\frac{1}{2}+\frac{\nu}{2}-}\big\|R^{-\frac{3}{2}}P_{\lambda_1}f P_{\lambda_2}g\big\|_{L^2}&\leq \sigma^{\frac{1}{2}+\frac{\nu}{2}-}\big\|P_{\lambda_1}f \big\|_{L^2}\big\|R^{-\frac{3}{2}}P_{\lambda_2}g\big\|_{L^\infty}\\&\lesssim \sigma^{\frac{1}{2}+\frac{\nu}{2}-}\big\|P_{\lambda_1}f \big\|_{L^2}\lambda_2\big\|P_{\lambda_2}g\big\|_{L^2}\\
&\lesssim \big(\frac{\sigma}{\lambda_1}\big)^{\frac{1}{2}+\frac{\nu}{2}-}\lambda_2^{-\frac{\nu}{2}+2\delta}\big\|P_{\lambda_1}f \big\|_{H_{\rho}^{\frac{1}{2}+\frac{\nu}{2}-}}\big\|P_{\lambda_2}g\big\|_{H_{\rho}^{1+\frac{\nu}{2}-2\delta-}}\\
\end{align*}
by Lemma 8.3 in \cite{KST0}. Again this is more than enough to square-sum over $\sigma, \lambda_1$ and sum over $\lambda_2$. These observations handle the case of small $\sigma$. We note that the $L^2$-type estimates for 
\[
R^{-j}fg,\,j\in \{0,1\}
\]
are just the same and in fact easier under the corresponding assumptions in the lemma.  
\\
Next, consider the case $\sigma\geq  \max\{\lambda_{1,2}\}$. If $\chi_{R<1}\nabla_R^l\big(R^{-\frac{3}{2}}P_{<\lambda}f\big)\in L^\infty$ uniformly in $\lambda>1$, then we get 
\[
\big\|\mathcal{L}^k\big(\chi_{R<1}R^{-\frac{3}{2}}P_{<\sigma}f P_{\lambda_2}g\big)\big\|_{L^2}\lesssim \lambda_2^k \|P_{\lambda_2}g\|_{L^2}
\]
Here we have used Lemma 8.4 in \cite{KST0}. It follows that 
\[
\big\|P_{\sigma}\big(\chi_{R<1}R^{-\frac{3}{2}}P_{<\sigma}f P_{\lambda_2}g\big)\big\|_{H_{\rho}^{\frac{1}{2}+\frac{\nu}{2}-}}\lesssim \big(\frac{\lambda_2}{\sigma}\big)^{k-\frac{1}{2}-\frac{\nu}{2}+} \|P_{\lambda_2}g\|_{H_{\rho}^{\frac{1}{2}+\frac{\nu}{2}-}}
\]
which suffices to square-sum over $\sigma$.
On the other hand, including a smooth cutoff $\chi_{R\geq 1}$, and assuming $\lambda_2\geq \lambda_1$ as we may, we get 
\begin{align*}
\big\|\mathcal{L}^k\big(\chi_{R\geq 1}R^{-\frac{3}{2}}P_{\lambda_1}f P_{\lambda_2}g\big)\big\|_{L^2}&\lesssim\sum_{m+l\leq k}\|\nabla_R^m P_{\lambda_1}f\|_{L^\infty} \| \nabla_R^l P_{\lambda_2}g\|_{L^2}\\
&\lesssim \sum_{m+l\leq k}\lambda_1^{-0+}\|\nabla_R^m P_{\lambda_1}f\|_{H_{\rho}^{\frac{1}{2}+\frac{\nu}{2}-}} \| \nabla_R^l P_{\lambda_2}g\|_{L^2}\\
&\lesssim \lambda_1^{-0+}\lambda_2^{-\frac{1}{2}-\frac{\nu}{2}+}\sum_{m+l\leq k}\lambda_1^m\lambda_2^l\|P_{\lambda_1}f\|_{H_{\rho}^{\frac{1}{2}+\frac{\nu}{2}-}} \|P_{\lambda_2}g\|_{H_{\rho}^{\frac{1}{2}+\frac{\nu}{2}-}}\\ 
\end{align*}
whence
\begin{align*}
\big\|P_{\sigma}\big(\chi_{R\geq 1}R^{-\frac{3}{2}}P_{\lambda_1}f P_{\lambda_2}g\big)\big\|_{H_{\rho}^{\frac{1}{2}+\frac{\nu}{2}-}}\lesssim \lambda_1^{-0+}\big(\frac{\lambda_2}{\sigma}\big)^{k-\frac{1}{2}-\frac{\nu}{2}+}\|P_{\lambda_1}f\|_{H_{\rho}^{\frac{1}{2}+\frac{\nu}{2}-}} \|P_{\lambda_2}g\|_{H_{\rho}^{\frac{1}{2}+\frac{\nu}{2}-}}.
\end{align*}
This again suffices to square-sum over $\sigma$ and $l^1$-sum over $\lambda_1$. 
If $g\in H^{1+\frac{\nu}{2}-2\delta-}$, we note that the argument for Lemma 8.5 in \cite{KST0} furnishes the bound 
\[
\big\|\mathcal{L}^k\big(R^{-\frac{3}{2}}P_{\lambda_1}f P_{\lambda_2}g\big)\big\|_{L^2}
\lesssim \lambda_1^{\frac{1}{2}}\lambda_2^k\|P_{\lambda_1}f\|_{H^{\frac{1}{2}}_{\rho}}\|P_{\lambda_2}g\|_{L_2^{\rho}},
\]
and so we get 
\[
\big\|\mathcal{L}^k P_{\sigma}\big(R^{-\frac{3}{2}}P_{\lambda_1}f P_{\lambda_2}g\big)\big\|_{H^{\frac{1}{2}+\frac{\nu}{2}-}_{\rho}}\lesssim \lambda_1^{\frac{1}{2}}\lambda_2^k\sigma^{\frac{1}{2}+\frac{\nu}{2}-}\|P_{\lambda_1}f\|_{H^{\frac{1}{2}}_{\rho}}\|P_{\lambda_2}g\|_{L_2^{\rho}}
\]
The duality argument in \cite{KST0} then yields (provided $\sigma>\lambda_2\geq \lambda_1$)
\[
\big\|P_{\sigma}\big(R^{-\frac{3}{2}}P_{\lambda_1}f P_{\lambda_2}g\big)\big\|_{H^{\frac{1}{2}+\frac{\nu}{2}-}_{\rho}}\lesssim
\big(\frac{\lambda_2}{\sigma}\big)^{\frac{1}{2}}\lambda_1^{-\frac{\nu}{2}+2\delta}\|P_{\lambda_1}f\|_{H^{\frac{1}{2}+\frac{\nu}{2}-}_{\rho}}\|P_{\lambda_2}g\|_{H^{1+\frac{\nu}{2}-2\delta-}_{\rho}}
\]
which suffices for the case $\lambda_1\leq \lambda_2<\sigma$, and the necessary summations. For the case $\lambda_1\geq \lambda_2$, one instead uses that 
\[
\big\|\mathcal{L}^k\big(R^{-\frac{3}{2}}P_{\lambda_1}f P_{\lambda_2}g\big)\big\|_{L^2}\lesssim \lambda_1^k\lambda_2\|P_{\lambda_1}f\|_{L^2}\|P_{\lambda_2}g\|_{L^2},
\]
which implies that 
\[
\big\|P_{\sigma}\big(R^{-\frac{3}{2}}P_{\lambda_1}f P_{\lambda_2}g\big)\big\|_{H_{\rho}^{\frac{1}{2}+\frac{\nu}{2}-}}\lesssim 
\big(\frac{\lambda_1}{\sigma}\big)^{k-\frac{1}{2}-\frac{\nu}{2}+} \lambda_2^{-\frac{\nu}{2}+2\delta}\|P_{\lambda_1}f\|_{H^{\frac{1}{2}+\frac{\nu}{2}-}_{\rho}}\|P_{\lambda_2}g\|_{H^{1+\frac{\nu}{2}-2\delta-}_{\rho}}
\]
This is again enough to sum over all dyadic frequencies. Finally, to obtain the inclusion $R^{-\frac{3}{2}}fg\in R^{\frac{3}{2}}L^\infty$, we observe that 
\begin{align*}
&|g(R)|=\\& \big|\int_0^\infty\phi(R,\xi)x(\xi)\rho(\xi)\,d\xi\big|\lesssim R^{\frac{3}{2}}\big(\int_0^\infty x^{2}(\xi)\langle \xi\rangle^{2+\nu-2\delta-}\rho(\xi)\,d\xi\big)^{\frac{1}{2}}\big(\int_0^\infty \langle \xi\rangle^{-2-\nu+2\delta+}\rho(\xi)\,d\xi\big)^{\frac{1}{2}}\\
& \lesssim \|g\|_{H^{1+\frac{\nu}{2}-2\delta-}_\rho}
\end{align*}
whence $|g(R)|\lesssim R^{\frac{3}{2}}\|g\|_{H^{1+\frac{\nu}{2}-2\delta-}_\rho}$. This implies 
\[
\|R^{-3}fg\|_{L^\infty}\lesssim \|f\|_{H_{\rho}^{\frac{1}{2}+\frac{\nu}{2}-}\cap R^{\frac{3}{2}}L^\infty}\|g\|_{H^{1+\frac{\nu}{2}-2\delta-}_\rho + R^{\frac{3}{2}}L^\infty}
\]
We also need to control $\|R^{-\frac{3}{2}}P_{<\lambda}\big(R^{-\frac{3}{2}}fg\big)\|_{L_x^\infty}$ for arbitrary dyadic $\lambda>1$. Write 
\begin{align}
&R^{-\frac{3}{2}}P_{<\lambda}\big(R^{-\frac{3}{2}}fg\big)\nonumber\\
&= R^{-\frac{3}{2}}P_{<\lambda}\big(\chi_{R\sim \tilde{R}}\tilde{R}^{-\frac{3}{2}}fg\big)\label{eq:RsimR_1}\\
&+R^{-\frac{3}{2}}P_{<\lambda}\big(\chi_{R\ll \tilde{R}}\tilde{R}^{-\frac{3}{2}}fg\big)\label{eq:RllR_1}\\
&+R^{-\frac{3}{2}}P_{<\lambda}\big(\chi_{R\gg\tilde{R}}\tilde{R}^{-\frac{3}{2}}fg\big)\label{eq:RggR_1}
\end{align}
for smooth cutoffs $\chi_{R\sim \tilde{R}}$ etc. To bound the first term on the right, we use that the operator $P_{<\lambda}$ is given by integration against the kernel 
\begin{equation}\label{eq:K<lambda}
K_{<\lambda}(R, \tilde{R}) = \chi_{R\sim \tilde{R}}\int_0^\infty \rho(\xi)\phi(R, \xi) \chi_{\xi<\lambda}\phi(\tilde{R}, \xi)\,d\xi
\end{equation}
for a smooth kernel function $ \chi_{\xi<\lambda}$. We claim that this kernel maps $L^\infty$ continuously into $L^\infty$. Taking this for granted, we obtain for the term \eqref{eq:RsimR_1} the bound 
\begin{align*}
\big\| R^{-\frac{3}{2}}P_{<\lambda}\big(\chi_{R\sim \tilde{R}}R^{-\frac{3}{2}}fg\big)\big\|_{L^\infty}
&\lesssim \sup_{\tilde{R}\sim 2^j}\big\|P_{<\lambda}\big(\chi_{\tilde{R}}\tilde{R}^{-3}fg\big)\big\|_{L^\infty}\\
&\lesssim\|f\|_{R^{\frac{3}{2}}L^\infty} \|g\|_{H_{\rho}^{1+\frac{\nu}{2}-2\delta-}+R^{\frac{3}{2}}L^\infty}
\end{align*}
 To get the $L^\infty$-boundedness of \eqref{eq:K<lambda}, write 
 \begin{align*}
&\chi_{R\sim \tilde{R}}\int_0^\infty \rho(\xi)\phi(R, \xi) \chi_{\xi<\lambda}\phi(\tilde{R}, \xi)\,d\xi\\
&= \sum_{N\,\text{dyadic}}\chi_{R\sim \tilde{R}\sim N}\int_0^\infty \rho(\xi)\phi(R, \xi) \chi_{\xi<\min\{\lambda, N^{-2}\}}\phi(\tilde{R}, \xi)\,d\xi\\
&+\sum_{N\,\text{dyadic}}\chi_{R\sim \tilde{R}\sim N}\int_0^\infty \rho(\xi)\phi(R, \xi) \chi_{ N^{-2}\leq \xi<\lambda}\phi(\tilde{R}, \xi)\,d\xi\\
\end{align*}
 Using Theorem~\ref{thm:spectral2}, one infers for the first term on the right the bound 
 \begin{align*}
 \big|\sum_{N\,\text{dyadic}}\chi_{R\sim \tilde{R}\sim N}\int_0^\infty \rho(\xi)\phi(R, \xi) \chi_{\xi<\min\{\lambda, N^{-2}\}}\phi(\tilde{R}, \xi)\,d\xi\big|\lesssim \frac{\chi_{R\sim \tilde{R}}}{R},
 \end{align*}
 and this kind of kernel is easily seen to act boundedly on $L^\infty$. For the oscillatory integral kernel above, write schematically, using 
Theorem~\ref{thm:spectral3}, Theorem~\ref{thm:spectral4}
\begin{align*}
&\chi_{R\sim \tilde{R}\sim N}\int_0^\infty \rho(\xi)\phi(R, \xi) \chi_{ N^{-2}\leq \xi<\lambda}\phi(\tilde{R}, \xi)\,d\xi\\
&= \chi_{R\sim \tilde{R}\sim N}\int_0^\infty \rho(\xi)a(\xi)^2\xi^{-\frac{1}{2}}e^{\pm i R\xi^{\frac{1}{2}}\pm i\tilde{R}\xi^{\frac{1}{2}}}\big(1+O(\frac{1}{R\xi^{\frac{1}{2}}})\big)^2\chi_{N^{-2}<\xi<\lambda}\,d\xi\\
& = \chi_{R\sim \tilde{R}\sim N}\big[- N\widehat{\chi_1}(N(\pm R\pm\tilde{R})) + \lambda\widehat{\chi_1}(\lambda(\pm R\pm\tilde{R}))\big] + O(\big|\log(\frac{R\pm \tilde{R}}{R})\big|\frac{\chi_{R\sim \tilde{R}\sim N}}{R}),
\end{align*}
for a suitable smooth and compactly supported function $\chi_1$, and the $L^\infty$-boundedness of the (sum over dyadic $N$ of) these operators follows easily. 
This concludes the estimate for \eqref{eq:RsimR_1}. 
\\
To bound the term \eqref{eq:RllR_1}, we break it into a number of constituents, using Theorem~\ref{thm:spectral1} - Theorem~\ref{thm:spectral4}. Write 
\begin{align*}
&R^{-\frac{3}{2}}P_{<\lambda}\big(\chi_{R\ll \tilde{R}}\tilde{R}^{-\frac{3}{2}}fg\big)\\
&= R^{-\frac{3}{2}}\int_0^\infty \int_0^\infty \chi_{R\ll \tilde{R}}\tilde{R}^{-\frac{3}{2}}f(\tilde{R})g(\tilde{R})\chi_{\xi<\lambda}\phi(\tilde{R}, \xi)\phi(R, \xi)\rho(\xi)\,d\xi d\tilde{R}
\end{align*}
with smooth cutoffs $\chi_{R\ll \tilde{R}}, \chi_{\xi<\lambda}$. We further split this as 
\begin{align}
&R^{-\frac{3}{2}}\int_0^\infty \int_0^\infty \chi_{R\ll \tilde{R}}\tilde{R}^{-\frac{3}{2}}f(\tilde{R})g(\tilde{R})\chi_{\xi<\lambda}\phi(\tilde{R}, \xi)\phi(R, \xi)\rho(\xi)\,d\xi d\tilde{R}\nonumber\\
&=R^{-\frac{3}{2}}\int_0^\infty \int_0^\infty \chi_{R\ll \tilde{R}}\chi_{R^2\xi\geq 1}\tilde{R}^{-\frac{3}{2}}f(\tilde{R})g(\tilde{R})\chi_{\xi<\lambda}\phi(\tilde{R}, \xi)\phi(R, \xi)\rho(\xi)\,d\xi d\tilde{R}\label{eq:tedious1}\\
&+R^{-\frac{3}{2}}\int_0^\infty \int_0^\infty \chi_{R\ll \tilde{R}}\chi_{R^{-2}>\xi\geq \tilde{R}^{-2}}\tilde{R}^{-\frac{3}{2}}f(\tilde{R})g(\tilde{R})\chi_{\xi<\lambda}\phi(\tilde{R}, \xi)\phi(R, \xi)\rho(\xi)\,d\xi d\tilde{R}\label{eq:tedious2}\\
&+R^{-\frac{3}{2}}\int_0^\infty \int_0^\infty \chi_{R\ll \tilde{R}}\chi_{\tilde{R}^2\xi<1}\tilde{R}^{-\frac{3}{2}}f(\tilde{R})g(\tilde{R})\chi_{\xi<\lambda}\phi(\tilde{R}, \xi)\phi(R, \xi)\rho(\xi)\,d\xi d\tilde{R}\label{eq:tedious3}
\end{align}
For the first term on the right, \eqref{eq:tedious1}, both functions $\phi(R, \xi)$, $\phi(\tilde{R}, \xi)$, are in the oscillatory regime, and can thus be written schematically as 
\[
\phi(R, \xi) = a(\xi)\xi^{-\frac{1}{4}}e^{\pm iR\xi^{\frac{1}{2}}}\big(1+O(\frac{1}{R\xi^{\frac{1}{2}}})\big),\,\phi(\tilde{R}, \xi) = a(\xi)\xi^{-\frac{1}{4}}e^{\pm i\tilde{R}\xi^{\frac{1}{2}}}\big(1+O(\frac{1}{\tilde{R}\xi^{\frac{1}{2}}})\big). 
\]
By applying integration by parts with respect to the variable $\xi^{\frac{1}{2}}$, we find 
\begin{align*}
\big|\eqref{eq:tedious1}\big|\lesssim  \int_0^\infty \chi_{R\ll \tilde{R}}(\frac{R}{\tilde{R}})^N \tilde{R}^{-4}\big|f(\tilde{R})\big|\big|g(\tilde{R})\big|\,d\tilde{R}
\end{align*}
and from here we get 
\begin{align*}
\big\|\eqref{eq:tedious1}\big\|_{L^\infty}\lesssim \big\|f\big\|_{R^{\frac{3}{2}}L^\infty}\big\|g\big\|_{H^{1+\frac{\nu}{2}-2\delta-}_{\rho}+R^{\frac{3}{2}}L^\infty}
\end{align*}
For the intermediate term \eqref{eq:tedious2}, one uses the expansions 
\[
\phi(R, \xi) = \phi_0(R) + \phi_0(R)O(R\xi^2),\,\phi(\tilde{R}, \xi) = a(\xi)\xi^{-\frac{1}{4}}e^{\pm i\tilde{R}\xi^{\frac{1}{2}}}\big(1+O(\frac{1}{\tilde{R}\xi^{\frac{1}{2}}})\big),
\]
and then uses again integration by parts with respect to $\xi^{\frac{1}{2}}$, obtaining bounds just as in the preceding case. 
Finally, for the remaining integral \eqref{eq:tedious3}, using the expansions 
\[
\phi(R, \xi) = \phi_0(R) + \phi_0(R)O(R\xi^2),\,\phi(\tilde{R}, \xi) = \phi_0(\tilde{R}) + \phi_0(\tilde{R})O(\tilde{R}\xi^2),
\]
we find 
\begin{align*}
\big|\eqref{eq:tedious3}\big|&\lesssim  \big(\int_0^\lambda \rho(\xi)\langle\xi\rangle^{-2}\,d\xi\big) \big\|\frac{f}{\tilde{R}^{\frac{3}{2}}}\big\|_{L^\infty}\big\|\frac{g}{\tilde{R}^{\frac{3}{2}}}\big\|_{L^\infty}\\
&\lesssim\log\lambda \big\|\frac{f}{\tilde{R}^{\frac{3}{2}}}\big\|_{L^\infty}\big\|\frac{g}{\tilde{R}^{\frac{3}{2}}}\big\|_{L^\infty}\\
\end{align*}
If we replace here the outer factor $R^{-\frac{3}{2}}$ by $R^{-1}$, one instead gets the bound 
\[
\lesssim  \big(\int_0^\lambda \rho(\xi)\langle\xi\rangle^{-\frac{9}{4}}\,d\xi\big) \big\|\frac{f}{\tilde{R}^{\frac{3}{2}}}\big\|_{L^\infty}\big\|\frac{g}{\tilde{R}^{\frac{3}{2}}}\big\|_{L^\infty},
\]
and so we no longer get a logarithmic correction for $\big\|R^{-1}P_{<\lambda}(R^{-\frac{3}{2}}fg)\big\|_{L^\infty}$. 

Observe that in order to bound $\|R^{-1}P_{<\lambda}R^{-1}fg\|_{L^\infty}$, and under the assumption $f\in RL^\infty, g\in RL^\infty$, proceeding just as before, we encounter instead of \eqref{eq:tedious3} a similar expression with the factors $R^{-\frac{3}{2}}$, $\tilde{R}^{-\frac{3}{2}}$ replaced by $R^{-1}, \tilde{R}^{-1}$. This we can then bound by 
\[
\lesssim \big\|\frac{f}{R}\big\|_{L^\infty}\big\|\frac{g}{R}\big\|_{L^\infty}\int_{R\ll \tilde{R}}\tilde{R}^{\frac{5}{2}}R^{\frac{1}{2}}\tilde{R}^{-4}\,d\tilde{R}\lesssim  \big\|\frac{f}{R}\big\|_{L^\infty}\big\|\frac{g}{R}\big\|_{L^\infty},
\]
thus without logarithmic correction. It is clear that the remaining cases occuring in the bound for \eqref{eq:RllR_1}, as well as for \eqref{eq:RsimR_1}, are easier for the expression $\big\|R^{-1}P_{<\lambda}(R^{-1}fg)\big\|_{L^\infty}$, and hence omitted. 
The bound for \eqref{eq:RggR_1} is more of the same. This completes the proof of the lemma. 
\end{proof}

\begin{lemma}\label{lem:multilin2}
Assume that all of $f, g, h$ are either in $H^{\frac{1}{2}+\frac{\nu}{2}-}_{\rho}\cap R^{\frac{3}{2}}L^\infty$ as well as with their frequency localized constituents $P_{<\lambda}(\cdot)\in \log \lambda R^{\frac{3}{2}}L^\infty$ and $\chi_{R<1}\nabla_R^l\big(R^{-\frac{3}{2}}P_{<\lambda}(\cdot)\big)\in L^\infty$, $l\geq 0$, uniformly in $\lambda>1$,  or in $H^{1+\frac{\nu}{2}-2\delta-}_{\rho}$. Then we have 
\[
R^{-3}fgh\in H^{\frac{1}{2}+\frac{\nu}{2}-2\delta-}\cap R^{\frac{3}{2}}L^\infty,\,P_{<\lambda}(R^{-3}fgh)\in \log\lambda R^{\frac{3}{2}}L^\infty,\,P_{<\lambda}(R^{-3}fgh)\in RL^\infty
\]
with the latter two inclusions uniformly in $\lambda>1$. Also, if $h_j\in H^{\frac{1}{2}+\frac{\nu}{2}-}_{\rho}\cap R^{\frac{3}{2}}L^\infty$ and further $P_{<\lambda}h_j\in RL^\infty$ as well as $\chi_{R<1}\nabla_R^l\big(R^{-1}P_{<\lambda}h_j\big)\in L^\infty$, $l\geq 0$, uniformly in $\lambda$, or else $h_j\in H_{\rho}^{1+\frac{\nu}{2}-2\delta-}$, for $j = 1,2,\ldots, 2N$, then we have 
\[
R^{-3}fgh \prod_{j=1}^N(\frac{1}{R}h_{2j} h_{2j-1})\in H^{\frac{1}{2}+\frac{\nu}{2}-2\delta-}
\]
We also get 
\[
R^{-\frac{3}{2}}fg\prod_{j=1}^N(\frac{1}{R}h_{2j} h_{2j-1})\in H^{\frac{1}{2}+\frac{\nu}{2}-\delta-}
\]
\end{lemma}
For the proof of this, one notes that by the preceding lemma, 
\[
R^{-\frac{3}{2}}fg \in  H^{\frac{1}{2}+\frac{\nu}{2}-\delta-}\cap R^{\frac{3}{2}}L^\infty,\,R^{-\frac{3}{2}}P_{<\lambda}(R^{-\frac{3}{2}}fg)\in \log\lambda L^\infty 
\]
uniformly in $\lambda>1$. Also, we have 
\[
R^{-1}P_{<\lambda}(R^{-\frac{3}{2}}fg)\in L^\infty 
\]
uniformly in $\lambda>1$. 
By another application of the preceding Lemma, we obtain the conclusions concerning $R^{-3}fgh$. The conclusion concerning 
\[
R^{-3}fgh \prod_{j=1}^N(\frac{1}{R}h_{2j} h_{2j-1})
\]
then follows by further iterative application of the preceding lemma. The last statement of the lemma follows similarly. 
\\
%To complete things, we need the following

%\begin{lemma}\label{lem:multilin3} Assume that all of $f\in H^{\frac{1}{2}+\frac{\nu}{2}}_{\rho}\cap R^{\frac{3}{2}}L^{\infty}\cap L^\infty$ and that $g, h$ are either in this space or else in $H^{1+\frac{\nu}{2}}_{\rho}$. Then we have 
%\[
%R^{-1}fgh\in H^{\frac{1}{2}+\frac{\nu}{2}}_{\rho}\cap R^{\frac{3}{2}}L^{\infty}\cap L^\infty
%\]
%\end{lemma}

%By comparison to the preceding lemma this offers new information only when $R\geq 1$, and is proved analogously to the preceding lemmas. 

We can now complete the estimate for the remaining two nonlinear source terms. Observe that we can write 
the first of these, \eqref{eq:N_2} in the form 
\begin{align*}
\frac{\sin(2u_{2k})}{2R^{\frac{3}{2}}}\big(1-\cos(2R^{-\frac{1}{2}}\tileps)\big) =& \frac{\sin(2u_{2k})}{2R^{\frac{3}{2}}}\big(R^{-\frac{1}{2}}\tileps\big)^2 q(R^{-1}\tileps^2)\\
&=\frac{\sin(2u_{2k})}{2R}R^{-\frac{3}{2}}\tileps^2q(R^{-1}\tileps^2)\\
\end{align*}
where $q(\cdot)$ is real analytic. 
By combining Lemma~\ref{lem:multilin2} and  Lemma~\ref{lem:8.1} (with $\alpha = \frac{1}{2}+\frac{\nu}{2}$) and using 
\[
\frac{\sin(2u_{2k})}{2R}\in IS(1, \mathcal{Q}),
\]
we find 

\begin{lemma}\label{lem:N_2Bound} We have the source term bound 
\[
\big\|\frac{\sin(2u_{2k})}{2R^{\frac{3}{2}}}\big(1-\cos(2R^{-\frac{1}{2}}\tileps)\big)\big\|_{H^{\frac{1}{2}+\frac{\nu}{2}-\delta-}_{\rho}}\lesssim \big\|\tileps\big\|_{H^{\frac{1}{2}+\frac{\nu}{2}-}_{\rho}\cap R^{\frac{3}{2}}L^\infty}^2
\]
provided we have 
\begin{equation}\label{eq:annoying1}
\big\|\tileps\big\|_{H^{\frac{1}{2}+\frac{\nu}{2}-}_{\rho}\cap R^{\frac{3}{2}}L^\infty}\lesssim 1,\,\big\|R^{-\frac{3}{2}}P_{<\lambda}\tileps\|_{L^\infty}\lesssim 1,\,\big\|\chi_{R<1}\nabla_R^l\big(R^{-\frac{3}{2}}P_{<\lambda}\tileps\big)\big\|_{L^\infty}\lesssim 1,\,l\geq 0
\end{equation}
uniformly in $\lambda>1$.
The same bound obtains with the space $H^{\frac{1}{2}+\frac{\nu}{2}-}_{\rho}\cap R^{\frac{3}{2}}L^\infty$ on the right replaced by $H^{1+\frac{\nu}{2}-2\delta-}_{\rho}$, and the bounds \eqref{eq:annoying1} replaced by 
\[
\|\tileps\|_{H^{1+\frac{\nu}{2}-2\delta-}_{\rho}}\lesssim 1.
\]
\end{lemma}

To deal with the last source term \eqref{eq:N_3}, we write 
\begin{align*}
\frac{\cos(2u_{2k})}{2R^{\frac{3}{2}}}\big(2R^{-\frac{1}{2}}\tileps - \sin(2R^{-\frac{1}{2}}\tileps)\big) = \cos(2u_{2k})\frac{\tileps^3}{R^3}q(R^{-1}\tileps^2)
\end{align*}
where again $q(\cdot)$ is real analytic. 
Combining Lemma~\ref{lem:multilin2}, and Lemma~\ref{lem:8.1}, we infer 
\begin{lemma}\label{lem:N_3Bound} We have the source term bound 
\[
\big\|\frac{\cos(2u_{2k})}{2R^{\frac{3}{2}}}\big(2R^{-\frac{1}{2}}\tileps - \sin(2R^{-\frac{1}{2}}\tileps)\big)\big\|_{H^{\frac{1}{2}+\frac{\nu}{2}-2\delta-}_{\rho}}\lesssim  \big\|\tileps\big\|_{H^{\frac{1}{2}+\frac{\nu}{2}-}_{\rho}\cap R^{\frac{3}{2}}L^\infty}^3
\]
provided we have 
\begin{equation}\label{eq:annoying2}
\big\|\tileps\big\|_{H^{\frac{1}{2}+\frac{\nu}{2}-}_{\rho}\cap R^{\frac{3}{2}}L^\infty}\lesssim 1,\,\big\|R^{-\frac{3}{2}}P_{<\lambda}\tileps\big\|_{L^\infty}\lesssim 1,\,\big\|\chi_{R<1}\nabla_R^l\big(R^{-\frac{3}{2}}P_{<\lambda}\tileps\big)\big\|_{L^\infty}\lesssim 1,\,l\geq 0
\end{equation}
uniformly in $\lambda>1$. The same bound obtains with the space $H^{\frac{1}{2}+\frac{\nu}{2}-}_{\rho}\cap R^{\frac{3}{2}}L^\infty$ on the right replaced by $H^{1+\frac{\nu}{2}-2\delta-}_{\rho}$, and \eqref{eq:annoying2} replaced by 
\[
\|\tileps\|_{H^{1+\frac{\nu}{2}-2\delta-}_{\rho}}\lesssim 1.
\]
\end{lemma}

\subsection{The first iterate} Recall that we have constructed the zeroth iterate via 
\[
x_0(\tau, \xi) = (U\lambda^{-2}\mathcal{F}\big[R^{\frac{1}{2}}\big(\tilde{e}_{2k-1}\big)\big])(\tau, \xi), 
\]
so that Proposition~\ref{prop:key} applies. 
Now we construct the {\it{first iterate}} via
\[
x_1(\tau, \xi) = (Uf_{0})(\tau, \xi),
\]
where we have 
\begin{align*}
-f_0 = &2\frac{\lambda_{\tau}}{\lambda}\mathcal{K}_0\big(\partial_{\tau} - \frac{\lambda_{\tau}}{\lambda}2\xi\partial_{\xi}\big)x_0 + (\frac{\lambda_{\tau}}{\lambda})^2\big[\mathcal{K}^2 - (\mathcal{K} - \mathcal{K}_0)^2 - 2[\xi\partial_{\xi}, \mathcal{K}_0]\big]x_0\\
&+\partial_{\tau}(\frac{\lambda_{\tau}}{\lambda})\mathcal{K}_0x_0+\lambda^{-2}\mathcal{F}\big[R^{\frac{1}{2}}\big(N_{2k-1}(R^{-\frac{1}{2}}\tilde{\eps}_0) + \tilde{e}_{2k-1}\big)\big] - c\tau^{-2}x_0
\end{align*}
Observe that we have 
\[
\big(\partial_{\tau} - \frac{\lambda_{\tau}}{\lambda}2\xi\partial_{\xi}\big)x_0\in \tau^{-N-1}L^{2, \frac{\nu}{2}-}_{\rho}
\]
Due to the smoothing property of $\mathcal{K}_0$, we conclude that 
\[
2\frac{\lambda_{\tau}}{\lambda}\mathcal{K}_0\big(\partial_{\tau} - \frac{\lambda_{\tau}}{\lambda}2\xi\partial_{\xi}\big)x_0\in \tau^{-N-2}L^{2, \frac{1}{2}+\frac{\nu}{2}-}_{\rho}
\]
Further, we get the even better bounds (which however we won't fully exploit) 
\[
 (\frac{\lambda_{\tau}}{\lambda})^2\big[\mathcal{K}^2 - (\mathcal{K} - \mathcal{K}_0)^2 - 2[\xi\partial_{\xi}, \mathcal{K}_0]\big]x_0\in \tau^{-N-2}L^{2, 1+\frac{\nu}{2}-}_{\rho}
\]
\[
\partial_{\tau}(\frac{\lambda_{\tau}}{\lambda})\mathcal{K}_0x_0 - c\tau^{-2}x_0\in \tau^{-N-2}L^{2, \frac{1}{2}+\frac{\nu}{2}-}_{\rho},
\]
while from Lemma~\ref{lem:N_2Bound}, Lemma~\ref{lem:N_3Bound} as well as \eqref{eq:N_1bound}, we infer 
\[
\big\|\lambda^{-2}\mathcal{F}\big[R^{\frac{1}{2}}\big(N_{2k-1}(R^{-\frac{1}{2}}\tilde{\eps}_0)\big\|_{L^{\frac{1}{2}+\frac{\nu}{2}-2\delta-}_{\rho}}\lesssim \tau^{-N-2}
\]
The key conclusion of all this is then the following 
\begin{lemma}\label{lem:Delta_1}The difference $\Delta x_1: = x_1 - x_0$ satisfies the bound 
\[
\big\|\Delta x_1(\tau, \cdot)\big\|_{L^{2,1+\frac{\nu}{2}-2\delta-}_{\rho}}\lesssim N^{-1}\tau^{-N},
\]
\[
\big\|\big(\partial_{\tau} - \frac{\lambda_{\tau}}{\lambda}2\xi\partial_{\xi}\big)\Delta x_1(\tau, \cdot)\big\|_{L^{2,\frac{1}{2}+\frac{\nu}{2}-2\delta-}_{\rho}}\lesssim N^{-1}\tau^{-N-1}
\]
The implicit constant is independent of $N$, whence picking $N$ large enough makes the overall constant on the right $\ll 1$. 
\end{lemma}

Note that the key aspect here is the gain of one derivative (which translates to a $1/2$ weight in terms of $\xi$). This is essential in order to replicate the reasoning used above for the new source term
\[
\lambda^{-2}\mathcal{F}\big[R^{\frac{1}{2}}\big(N_{2k-1}(R^{-\frac{1}{2}}\tilde{\eps}_1)\big]
\]
where we define the first iterate on the physical side via 
\begin{align*}
\tileps_1(\tau, R) = \int_0^\infty \phi(R, \xi)x_1(\tau, \xi)\rho(\xi)\,d\xi=&\int_0^\infty \phi(R, \xi)\Delta x_1(\tau, \xi)\rho(\xi)\,d\xi\\
&+\int_0^\infty \phi(R, \xi)x_0(\tau, \xi)\rho(\xi)\,d\xi\\
\end{align*}
Thus from Proposition~\ref{prop:key}, the remark following it, as well as Corollary~\ref{cor:refine} and the preceding lemma, we infer that we can write 
\[
\tileps_1(\tau, \cdot) = \tileps_1^{(1)}(\tau, \cdot) + \tileps_1^{(2)}(\tau, \cdot),
\]
where we have 
\[
\tileps_1^{(1)}(\tau, \cdot)\in \tau^{-N}\big(H^{\frac{1}{2}+\frac{\nu}{2}-}_{\rho}\cap R^{\frac{3}{2}}L^{\infty}\big), \nabla_R^l\big(R^{-\frac{3}{2}}P_{<\lambda}\tileps_1^{(1)}(\tau, \cdot)\big)\in \tau^{-N}L^\infty,\,l\geq 0,
\]
the latter inclusion uniformly in $\lambda>1$, while we have 
\[
\tileps_1^{(2)}(\tau, \cdot)\in \tau^{-N}H^{1+\frac{\nu}{2}-2\delta-}_{\rho}
\]
This is precisely the kind of structure necessary to invoke the bound \eqref{eq:N_1bound} as well as Lemma~\ref{lem:N_2Bound}, Lemma~\ref{lem:N_3Bound}. 

\subsection{Higher iterates}

Here we have 
\[
x_j(\tau, \xi) = (Uf_{j-1})(\tau, \xi),\,j\geq 2,
\]
and we have 
\begin{align*}
-f_{j-1} = &2\frac{\lambda_{\tau}}{\lambda}\mathcal{K}_0\big(\partial_{\tau} - \frac{\lambda_{\tau}}{\lambda}2\xi\partial_{\xi}\big)x_{j-1} + (\frac{\lambda_{\tau}}{\lambda})^2\big[\mathcal{K}^2 - (\mathcal{K} - \mathcal{K}_0)^2 - 2[\xi\partial_{\xi}, \mathcal{K}_0]\big]x_{j-1}\\
&+\partial_{\tau}(\frac{\lambda_{\tau}}{\lambda})\mathcal{K}_0x_{j-1}+\lambda^{-2}\mathcal{F}\big[R^{\frac{1}{2}}\big(N_{2k-1}(R^{-\frac{1}{2}}\tilde{\eps}_{j-1}) + \tilde{e}_{2k-1}\big)\big] - c\tau^{-2}x_{j-1}
\end{align*}
Then using induction on $j$ and exactly the same bounds as in the preceding subsection, one infers with 
\[
\Delta x_j = x_j - x_{j-1}
\]
the bounds 
\[
\big\|\Delta x_j(\tau, \cdot)\big\|_{L^{2,1+\frac{\nu}{2}-2\delta-}_{\rho}}\lesssim N^{-j}\tau^{-N},
\]
\[
\big\|\big(\partial_{\tau} - \frac{\lambda_{\tau}}{\lambda}2\xi\partial_{\xi}\big)\Delta x_j(\tau, \cdot)\big\|_{L^{2,\frac{1}{2}+\frac{\nu}{2}-2\delta-}_{\rho}}\lesssim N^{-j}\tau^{-N-1}
\]
The desired fixed point of \eqref{eq:Fourier3} is now obtained via 
\[
x(\tau, \xi) = x_0(\tau, \xi) + \sum_{j=1}^\infty \Delta x_j(\tau, \xi)
\]
and is a function in $H^{\frac{1}{2}+\frac{\nu}{2}-}_{\rho}$, such that $\partial_{\tau}x(\tau, \cdot)\in H^{\frac{\nu}{2}-}_{\rho}$. 
Due to Lemma 9.1 of \cite{KST0}, the corresponding
\[
\eps(\tau, R): = R^{-\frac{1}{2}}\int_0^\infty \phi(R, \xi)x(\tau, \xi)\rho(\xi)\,d\xi
\]
satisfies $\eps(\tau, \cdot)\in \tau^{-N}H^{1+\nu-}_{\R^2}$,  as well as $\partial_{\tau}\eps(\tau, \cdot)\in \tau^{-N-1}H^{\nu-}_{\R^2}$. 
This is the desired solution.

\bigskip

% Enter the first author's name and address:

\centerline{\scshape Can Gao }
\medskip
{\footnotesize
% please put the address of the first author
 \centerline{B\^{a}timent des Math\'ematiques, EPFL}
\centerline{Station 8, 
CH-1015 Lausanne, 
  Switzerland}
  \centerline{\email{can.gao@epfl.ch}}
} % Do not forget to end the {\footnotesize by the sign }

\medskip

\centerline{\scshape Joachim Krieger }
\medskip
{\footnotesize
% please put the address of the first author
 \centerline{B\^{a}timent des Math\'ematiques, EPFL}
\centerline{Station 8, 
CH-1015 Lausanne, 
  Switzerland}
  \centerline{\email{joachim.krieger@epfl.ch}}
} % Do not forget to end the {\footnotesize by the sign }

\end{document}